\documentclass[11pt]{article}
\usepackage{amssymb,amsmath,amsfonts}
\usepackage{graphicx}
\usepackage{color}
\usepackage[mathscr]{eucal}

\allowdisplaybreaks  

\newcommand{\sinush}{\mathop{\rm sinh}}
\newcommand{\cosinh}{\mathop{\rm cosh}}

\newcommand{\arctanh}{\mathop{\rm artanh}}

\newcommand{\tang}{\mathop{\rm tan}}

\newtheorem{lemma}{Lemma}
\newtheorem{theorem}{Theorem}

\newtheorem{remark}{Remark}

\begin{document}

\title{\large{A NECESSARY FLEXIBILITY CONDITION OF A NONDEGENERATE SUSPENSION IN LOBACHEVSKY 3-SPACE}}

\author{Dmitriy~Slutskiy \thanks{The author is supported in part by the Council of
Grants from the President of the Russian Federation (Grant
NSh-6613.2010.1), by the Federal Targeted Programme on Scientific
and Pedagogical-Scientific Staff of Innovative Russia for 2009-2013
(State Contract No. 02.740.11.0457) and by the Russian Foundation
for Basic Research (Grant 10-01-91000-ANF\_a).}}

\date{}

\maketitle


\begin{abstract}

We show that some combination of the lengths of all edges of the
equator of a flexible suspension in Lobachevsky 3-space is equal to
zero (each length is taken either positive or negative in this
combination).

\end{abstract}

\begin{center}

\end{center}

\section{Introduction} \label{introduction}

A polyhedron (more precisely, a polyhedral surface) is said to be
flexible if its spatial shape can be changed continuously due to
changes of its dihedral angles only, i.e., if every face remains
congruent to itself during the flex.

In 1897 R.~Bricard \cite{Bricard1897} described all flexible
octahedra in Euclidean 3-space. The Bricard's octahedra were the
first examples of flexible polyhedra (with
self-in\-ter\-sec\-ti\-ons). BricardTs octahedra are special cases
of Euclidean flexible suspensions. In 1974 R.~Connelly
\cite{Connelly1980} proved that some combination of the lengths of
all edges of the equator of a flexible suspension in Euclidean
3-space is equal to zero (each length is taken either positive or
negative in this combination). The method applied by R.~Connelly, is
to reduce the problem to the study of an analytic function of
complex variable in neighborhoods of its singular points.

In 2001 S.\,N.~Mikhalev \cite{Mikhalev2001} reproved the
above-mentioned result of R.~Connelly by algebraic methods.
Moreover, S.\,N.~Mikhalev proved that for every spatial
quadrilateral formed by edges of a flexible suspension and
containing its both poles there is a combination of the lengths
(taken either positive or negative) of the edges of the
quadrilateral, which is equal to zero.

The aim of this work is to prove a similar result for the equator of
a flexible suspension in Lobachevsky 3-space, applying the method of
Connelly \cite{Connelly1980}.

\section{Formulating the main result} \label{The_main_result's_formulation}

Let $\mathscr{K}$ be a simplicial complex. A \emph{polyhedron} (a
\emph{polyhedral surface}) in Lo\-ba\-chevs\-ky 3-space is a
continuous map from $\mathscr{K}$ to $\mathbb{H}^3$, which sends
every $k$-dimensional simplex of $\mathscr{K}$ into a subset of a
$k$-dimensional plane of Lobachevsky space $(k\leq2)$. Images of
topological $2$-simplices are called faces, images of topological
$1$-simplices are called edges and images of topological
$0$-simplices are called vertices of the polyhedron. Note that in
our definition an image of a simplex can be degenerate (for
instance, a face can lie on a straight hyperbolic line, and an edge
can be reduced to one point), and faces can intersect in their
interior points. If $v_1,...,v_W$ are the vertices of $\mathscr{K}$,
and if $\mathscr{P}:\mathscr{K}\rightarrow\mathbb{H}^3$ is a
polyhedron, then $\mathscr{P}$ is determined by $W$ points
$P_1,...,P_W\in\mathbb{H}^3$, where
$P_j\stackrel{\mathrm{def}}{=}\mathscr{P}(v_j)$, $j=1,...,W$.

If $\mathscr{P}:\mathscr{K}\rightarrow\mathbb{H}^3$ and
$\mathscr{Q}:\mathscr{K}\rightarrow\mathbb{H}^3$ are two polyhedra,
then we say $\mathscr{P}$ and $\mathscr{Q}$ are \emph{congruent} if
there exists a motion
$\mathscr{A}:\mathbb{H}^3\rightarrow\mathbb{H}^3$ such that
$\mathscr{Q}=\mathscr{A}\circ\mathscr{P}$ (i.e. the isometric
mapping $\mathscr{A}$ sends every vertex of $\mathscr{P}$ into a
corresponding vertex of $\mathscr{Q}$: $Q_j=\mathscr{A}(P_j)$, or in
other words $\mathscr{Q}(v_j)=\mathscr{A}(\mathscr{P}(v_j))$, $j =
1,...,W$). We say $\mathscr{P}$ and $\mathscr{Q}$ are
\emph{isometric} (\emph{in the intrinsic metric}) if each edge of
$\mathscr{P}$ has the same length as the corresponding edge of
$\mathscr{Q}$, i.e. if $\langle v_j,v_k\rangle$ is a $1$-simplex of
$\mathscr{K}$ then
$\mathrm{d}_{\mathbb{H}^3}(Q_j,Q_k)=\mathrm{d}_{\mathbb{H}^3}(P_j,P_k)$,
where $\mathrm{d}_{\mathbb{H}^3}(\cdot,\cdot)$ stands for the
distance in Lobachevsky space $\mathbb{H}^3$.

A polyhedron $\mathscr{P}$ is \emph{flexible} if, for some
continuous one parameter family of polyhedra
$\mathscr{P}_t:\mathscr{K}\rightarrow\mathbb{H}^3$, $0\leq t\leq1$,
the following three conditions hold true: $(1)$
$\mathscr{P}_0=\mathscr{P}$; $(2)$ each $\mathscr{P}_t$ is isometric
to $\mathscr{P}_0$; $(3)$ some $\mathscr{P}_t$ is not congruent to
$\mathscr{P}_0$.

Let $\mathscr{K}$ be defined as follows: $\mathscr{K}$ has vertices
$v_0,v_1,...,v_V,v_{V+1}$, where $v_1,...,v_V$ form a cycle ($v_j$
adjacent to $v_{j+1}$, $j=1,...,V-1$, and $v_V$ adjacent to $v_1$),
and $v_0$ and $v_{V+1}$ are each adjacent to all of $v_1,...,v_V$.
Each polyhedron $\mathscr{P}$ based on $\mathscr{K}$ is called a
\emph{suspension}. Call
$N\stackrel{\mathrm{def}}{=}\mathscr{P}(v_0)$ the north pole, and
$S\stackrel{\mathrm{def}}{=}\mathscr{P}(v_{V+1})$ the south pole,
and $P_j\stackrel{\mathrm{def}}{=}\mathscr{P}(v_j)$, $j=1,...,V$
vertices of the \emph{equator} $\mathscr{P}$.

\begin{figure}
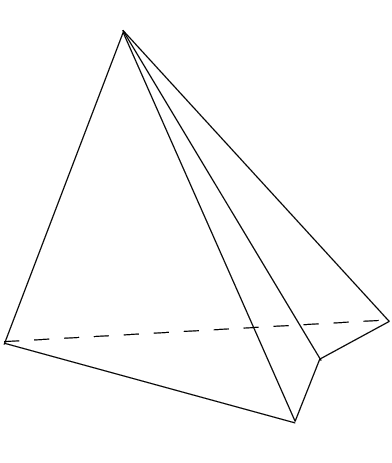 \hfill
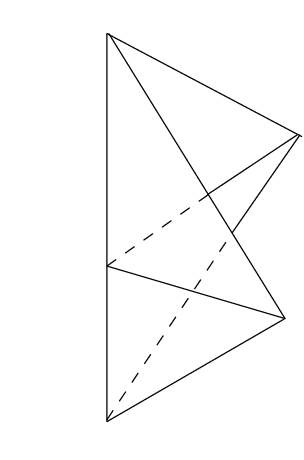 \\
\parbox[t]{7.5cm}{\caption{A double covered cap.}\label{DoubleCoveredCap}} \hfill
\parbox[t]{7cm}{\caption{A suspension with a wing.}\label{SuspensionWithWings}}
\end{figure}

Assume that a suspension $\mathscr{P}$ is flexible. If we suppose
the segment $NS$ to be an extra edge, then $\mathscr{P}$ becomes a
set of $V$ tetrahedra glued cyclically along their common edge $NS$.
We call a suspension \emph{nondegenerate} if none of these
tetrahedra lies on a hyperbolic $2$-plane. Note that a nondegenerate
suspension $\mathscr{P}$ does not flex if the distance between $N$
and $S$ remains constant. Therefore, as in the Euclidean case
\cite{Connelly1980} we assume that the length of $NS$ is variable
during the flex of $\mathscr{P}$. Examples of degenerate suspensions
are a double covered cap
--- a suspension with coinciding poles
(see Fig.~\ref{DoubleCoveredCap}), and a suspension with a wing
--- a suspension whose vertices $N$, $S$, $P_{i-1}$, and $P_{i+1}$
lie on a straight line for some $i$ (see
Fig.~\ref{SuspensionWithWings}). In this paper we will not study the
degenerate flexible suspensions.

The main result of the paper is
\begin{theorem}\label{sum_of_equator_edges_theorem}
Let $\mathscr{P}$ be a nondegenerate flexible suspension in
Lobachevsky \linebreak $3$-space with the poles $S$ and $N$, and
with the vertices of the equator $P_{j}$, $j=1,...,V$. Then for some
set of signs $\sigma_{j,j+1}\in\{+1,-1\}$, $j=1,...,V$, the
combination of the lengths $e_{j,j+1}$ of all edges $P_{j}P_{j+1}$
of the equator of $\mathscr{P}$ taken with the corresponding signs
$\sigma_{j,j+1}$ is equal to zero, i.e.
\begin{equation}\label{flex_eqn_at_infty_3}
\sum\limits_{j=1}^{V}\sigma_{j,j+1}e_{j,j+1}=0.
\end{equation}
\emph{(}Here and below, by definition, it is considered that
$P_{V+1}\stackrel{\mathrm{def}}{=}P_{1}$,
$P_{V}P_{V+1}\stackrel{\mathrm{def}}{=}P_{V}P_{1}$,
$\sigma_{V,V+1}\stackrel{\mathrm{def}}{=}\sigma_{V,1}$, and
$e_{V,V+1}\stackrel{\mathrm{def}}{=}e_{V,1}$.\emph{)}
\end{theorem}

\section{Connelly's equation of flexibility of a suspension}
\label{Connelly's_equation_on_the_suspension_flexibility}

R.~Connelly in \cite{Connelly1980} obtained an equation of
flexibility of a nondegenerate suspension in Euclidean 3-space.
Following him, in this section we will obtain an equation of
flexibility of a nondegenerate suspension in Lobachevsky 3-space.

\begin{figure}
\includegraphics{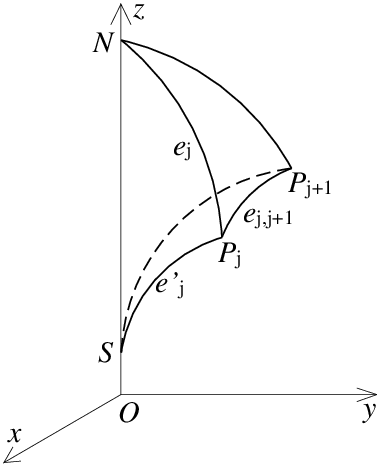} \hfill
\includegraphics{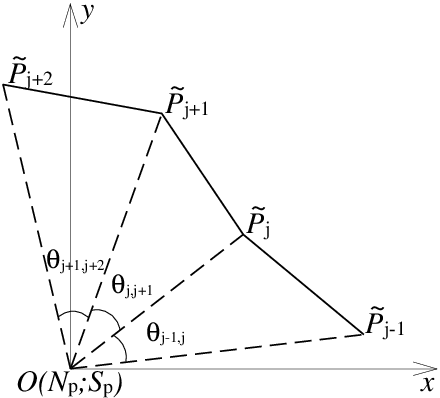} \\
\parbox[t]{7.5cm}{ \caption{A fragment of the lateral surface of $\mathscr{P}$.}\label{SuspensionPartInH3}} \hfill
\parbox[t]{7cm}{\caption{A projection of $\mathscr{P}$ on $Oxy$ .}\label{SuspensionProjection}}
\end{figure}

Let us place a nondegenerate suspension $\mathscr{P}$ into the
Poincar\'e upper half-space model \cite{Anderson2005} of Lobachevsky
3-space $\mathbb{H}^3$ in such a way that the poles $N$ and $S$ of
$\mathscr{P}$ lie on the axis $Oz$ of the Cartesian coordinate
system of the Poincar\'e model (see Fig.~\ref{SuspensionPartInH3}).
Let $S$ has the coordinates $(0,0,z_S)$, $N$ has the coordinates
$(0,0,z_N)$, and $P_j$ has the coordinates $(x_j,y_j,z_j)$,
$j=1,...,V$. Also we denote the length of the edge $NP_{j}$ by
$e_{j}$, and the length of $SP_{j}$ by $e_{j}'$, $j=1,...,V$.

Consider a Euclidean orthogonal projection $\widetilde{\mathscr{P}}$
of $\mathscr{P}$ on the plane $Oxy$ (see
Fig.~\ref{SuspensionProjection}). Also $\widetilde{\mathscr{P}}$ is
a hyperbolic projection of $\mathscr{P}$ on $Oxy$ from the only
point at infinity of $\mathbb{H}^3$ which does not lie on $Oxy$.
This projection sends poles $N$ and $S$ of $\mathscr{P}$ to the
origin $O$ $(0,0)$ on the plane $Oxy$, $P_j$ to the point
$\widetilde{P}_j$ $(x_j,y_j)$, edges $NP_j$ and $SP_j$ to the
Euclidean segment $O\widetilde{P}_j$, and the egde $P_{j}P_{j+1}$ of
the equator of $\mathscr{P}$ to the Euclidean segment
$\widetilde{P}_{j}\widetilde{P}_{j+1}$, $j=1,...,V$ (here and below
$\widetilde{P}_{V+1}\stackrel{\mathrm{def}}{=}\widetilde{P}_{1}$,
$x_{V+1}\stackrel{\mathrm{def}}{=}x_{1}$,
$y_{V+1}\stackrel{\mathrm{def}}{=}y_{1}$,
$z_{V+1}\stackrel{\mathrm{def}}{=}z_{1}$).

\begin{figure}
\begin{center}
\includegraphics{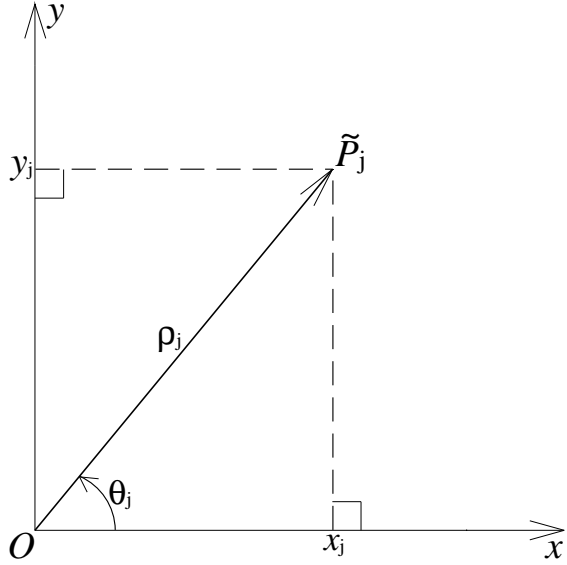}
\caption{The coordinates of
$\widetilde{P}_{j}$.}\label{PointCoordinates}
\end{center}
\end{figure}

Polar coordinates $(\rho_j,\theta_j)$ of $\widetilde{P}_j$,
$j=1,...,V$, are related to its Cartesian coordinates by the
formulas (see Fig.~\ref{PointCoordinates}):
\begin{equation}\label{rho_j_theta_j}
\rho_{j}=\sqrt{x_{j}^2+y_{j}^2},\quad
\sin\theta_{j}=\frac{y_j}{\rho_{j}}=\frac{y_j}{\sqrt{x_{j}^2+y_{j}^2}},\quad
\cos\theta_{j}=\frac{x_j}{\rho_{j}}=\frac{x_j}{\sqrt{x_{j}^2+y_{j}^2}}.
\end{equation}

Note that by construction, the dihedral angle $\theta_{j,j+1}$ of
the tetrahedron $NSP_{j}P_{j+1}$ at the edge $NS$ is equal to the
flat angle $\angle\widetilde{P}_{j}O\widetilde{P}_{j+1}$,
$j=1,...,V$, and
\begin{equation}\label{theta_j_j+1}
\theta_{j,j+1}=\theta_{j+1}-\theta_{j}.
\end{equation}
Note as well that the value of $\theta_{j,j+1}$ can be negative.
Applying the trigonometric ratio of the difference of two angles
and~(\ref{theta_j_j+1}), we get:
\begin{equation}\label{cos_theta_j_j+1}
\cos\theta_{j,j+1}=\cos\theta_{j+1}\cos\theta_{j}
+\sin\theta_{j+1}\sin\theta_{j},\quad
\sin\theta_{j,j+1}=\sin\theta_{j+1}\cos\theta_{j}
-\cos\theta_{j+1}\sin\theta_{j}.
\end{equation}
Taking into account (\ref{rho_j_theta_j}) we reduce
(\ref{cos_theta_j_j+1}) to
\begin{equation*}\label{cos_theta_j_j+1_by_x_y}
\cos\theta_{j,j+1}
=\frac{x_{j}x_{j+1}+y_{j}y_{j+1}}{\sqrt{x_{j+1}^2+y_{j+1}^2}\sqrt{x_{j}^2+y_{j}^2}},\quad
\sin\theta_{j,j+1}=
\frac{x_{j}y_{j+1}-y_{j}x_{j+1}}{\sqrt{x_{j+1}^2+y_{j+1}^2}\sqrt{x_{j}^2+y_{j}^2}}.
\end{equation*}

Then, according to Euler's formula,
\begin{equation}\label{exp(i_theta)}
e^{i\theta_{j,j+1}}=\cos\theta_{j,j+1}+i\sin\theta_{j,j+1}=
\frac{(x_{j}x_{j+1}+y_{j}y_{j+1})+i(x_{j}y_{j+1}-y_{j}x_{j+1})}
{\sqrt{x_{j+1}^2+y_{j+1}^2}\sqrt{x_{j}^2+y_{j}^2}}.
\end{equation}

Following R.~Connelly \cite{Connelly1980}, we remark that the sum of
the dihedral angles $\theta_{j,j+1}$ of all tetrahedra
$NSP_{j}P_{j+1}$, $j=1,...,V$, at the edge $NS$ is constant and a
multiple of $2\pi$ (here and below
$\theta_{V,V+1}\stackrel{\mathrm{def}}{=}\theta_{V,1}$,
$\theta_{V+1}\stackrel{\mathrm{def}}{=}\theta_{1}$,
$\rho_{V+1}\stackrel{\mathrm{def}}{=}\rho_{1}$), i.e.
\begin{equation}\label{full angle_thetas}
\sum_{j=1}^{V}\theta_{j,j+1}=2\pi m \quad \mbox{for some integer}
\,\,\, m,
\end{equation}
and remains so during the deformation of the suspension, when the
values of the angles $\theta_{j,j+1}$, $j=1,...,V$, vary
continuously.

We rewrite the equation of flexibility~(\ref{full angle_thetas}) in
a convenient form:
\begin{equation}\label{full angle_exp}
\prod_{j=1}^{V}e^{i\theta_{j,j+1}}=1.
\end{equation}
Thus, taking into account~(\ref{exp(i_theta)}), we see that
coordinates of vertices of $\mathscr{P}$ are related as follows:
\begin{equation}\label{full angle_coords}
\prod_{j=1}^{V}\frac{(x_{j}x_{j+1}+y_{j}y_{j+1})+i(x_{j}y_{j+1}-y_{j}x_{j+1})}
{x_{j}^2+y_{j}^2}=1,
\end{equation}
or in other notations
\begin{equation}\label{full angle_functions}
\prod_{j=1}^{V}F_{j,j+1}=\prod_{j=1}^{V}\frac{G_{j,j+1}}
{\rho_{j}\rho_{j+1}}=\prod_{j=1}^{V}\frac{G_{j,j+1}} {\rho_{j}^2}=1,
\end{equation}
where
$G_{j,m}=(x_{j}x_{m}+y_{j}y_{m})+i(x_{j}y_{m}-y_{j}x_{m})$,$F_{j,m}=\frac{G_{j,m}}
{\rho_{j}\rho_{m}}$, $j,m=1,...,V$, and
$G_{V,V+1}\stackrel{\mathrm{def}}{=}G_{V,1}$,
$F_{V,V+1}\stackrel{\mathrm{def}}{=}F_{V,1}$.

When studying the deformation $\mathscr{P}_t$ of the suspension
$\mathscr{P}$, all objects and values related to $\mathscr{P}_t$
naturally succeed from the notations for the corresponding entities
related to $\mathscr{P}$. For example, the coordinate $x_{j}(t)$ of
the point $P_{j}(t)$ of the deformation $\mathscr{P}_t$ corresponds
to the coordinate $x_{j}$ of the point $P_{j}$ of the suspension
$\mathscr{P}$, the dihedral angle $\theta_{j,j+1}(t)$ of the
tetrahedron $N(t)S(t)P_{j}(t)P_{j+1}(t)$ at the edge $N(t)S(t)$
corresponds to the dihedral angle $\theta_{j,j+1}$ of the
tetrahedron $NSP_{j}P_{j+1}$ at the edge $NS$, etc.

\section{The equation of flexibility of a suspension in terms of the lengths of its edges}
\label{Dependence_of_the_flexibility_equation_of_a_suspension_on_the_lengths_of_its_edges}

\begin{figure}
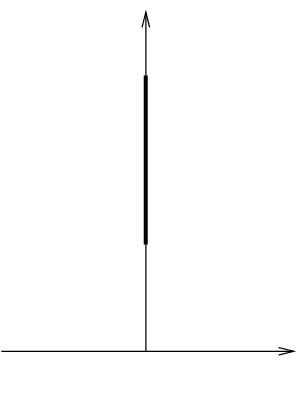 \hfill
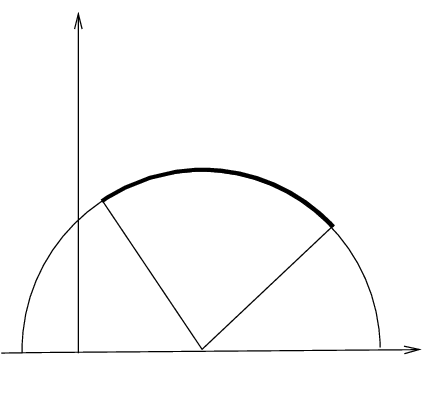 \\
\parbox[t]{7cm}{ \caption{Points on a plane in the lemma~\ref{vert_points_distance_lemma}.}\label{fig_lemma1}} \hfill
\parbox[t]{7cm}{\caption{Points on a plane in the lemma~\ref{points_distance_lemma}.}\label{fig_lemma2}}
\end{figure}

In this section we are going to express the equation of flexibility
of a suspension~(\ref{full angle_coords}) in terms of the lengths of
edges of $\mathscr{P}$. Recall that the lengths of the edges of
$\mathscr{P}$ remain constant during the flex. To this purpose we
need to demonstrate the truth of two following statements. The first
of them can be verified by direct calculation (see also
Fig.~\ref{fig_lemma1}).

\begin{lemma}\label{vert_points_distance_lemma}
Given a Poincar\'e upper half-plane $\mathbb{H}^2$ with the
coordinates $(\rho,z)$ \emph{(}i.e., with the metric given by the
formula $ds^2=\frac{d\rho^2+dz^2}{z^2}$\emph{)}. Then the distance
between the points $A$ $(\rho_{0},z_{A})$ and $B$
$(\rho_{0},z_{B})$, having the same first coordinate $\rho_{0}$, is
calculated by the formula
\begin{equation}\label{vert_points_distance}
\mathrm{d}_{\mathbb{H}^2}(A,B)=\Big|\ln\frac{z_{B}}{z_{A}}\Big|.
\end{equation}
\end{lemma}

\begin{lemma}\label{points_distance_lemma}
Given a Poincar\'e upper half-plane $\mathbb{H}^2$ with the
coordinates $(\rho,z)$ \emph{(}i.e., with the metric given by the
formula $ds^2=\frac{d\rho^2+dz^2}{z^2}$\emph{)}. Then the distance
$l\stackrel{\mathrm{def}}{=}\mathrm{d}_{\mathbb{H}^2}(A,B)$ between
the points $A$ $(\rho_{A},z_{A})$ and $B$ $(\rho_{B},z_{B})$ is
related to their coordinates by the formula
\begin{equation}\label{points_distance_coords}
(\rho_{B}-\rho_{A})^2+z_{A}^2+z_{B}^2=2z_{A}z_{B}\cosinh l.
\end{equation}
\end{lemma}

\textbf{Proof.} According to the part~$(2)$ of the
Corollary~$A.5.8$~\cite{BP2003}, the distance between the points
with the coordinates $(x,t)$ and $(y,s)$ in the Poincar\'e upper
half-space model $\mathbb{R}^{n}\times\mathbb{R}^{+}$ of Lobachevsky
$(n+1)$-space $\mathbb{H}^{n+1}$ is calculated by the formula
\begin{equation}\label{distance_in_H^n+1_BP}
\mathrm{d}_{\mathbb{H}^{n+1}}((x,t),(y,s))=2\arctanh\bigg(
\frac{\|x-y\|^2+(t-s)^2}{\|x-y\|^2+(t+s)^2}\bigg)^{1/2},
\end{equation}
where the symbol $\|\cdot\|$ stands for the standard Euclidean norm
in $\mathbb{R}^{n}$.

By~(\ref{distance_in_H^n+1_BP}) the distance between the points $A$
and $B$ (see Fig.~\ref{fig_lemma2}) is calculated by the formula
\begin{equation}\label{d_H^2(A,B)_BP}
l=2\arctanh\bigg(
\frac{(\rho_{A}-\rho_{B})^2+(z_{A}-z_{B})^2}{(\rho_{A}-\rho_{B})^2+(z_{A}+z_{B})^2}\bigg)^{1/2},
\end{equation}
where $n=1$, $(x,t)=(\rho_{A},z_{A})$ and $(y,s)=(\rho_{B},z_{B})$.

After a series of transformations of the
formula~(\ref{d_H^2(A,B)_BP}) we get:
\begin{equation}\label{d_H^2(A,B)_BP_3}
(\rho_{A}-\rho_{B})^2\Big({\cosinh}^{2}\frac{l}{2}-{\sinush}^{2}\frac{l}{2}\Big)+
(z_{A}^2+z_{B}^2)\Big({\cosinh}^{2}\frac{l}{2}-{\sinush}^{2}\frac{l}{2}\Big)=
2z_{A}z_{B}\Big({\cosinh}^{2}\frac{l}{2}+{\sinush}^{2}\frac{l}{2}\Big).
\end{equation}
By two identities of hyperbolic geometry,
${\cosinh}^{2}\frac{l}{2}-{\sinush}^{2}\frac{l}{2}=1$ and $\cosinh
l={\cosinh}^{2}\frac{l}{2}+{\sinush}^{2}\frac{l}{2}$,
(\ref{d_H^2(A,B)_BP_3}) reduces to~(\ref{points_distance_coords}).
$\square$

Let us express $G_{j,j+1}$ and $\rho_{j}^2$ in terms of the length
of edges of $\mathscr{P}$.

We assume that the coordinates of the south pole $S$ are $(0,0,1)$.
Let $t\stackrel{\mathrm{def}}{=}e^{\mathrm{d}_{\mathbb{H}^3}(N,S)}$,
where $\mathrm{d}_{\mathbb{H}^3}(N,S)$ is the distance between the
poles $N$ and $S$ of $\mathscr{P}$. Without loss of generality, we
assume that $z_{N}\geq z_{S}$. Then, by
Lemma~\ref{vert_points_distance_lemma}, the coordinates of $N$ are
$(0,0,t)$.

Applying Lemma~\ref{points_distance_lemma} to the points $S$ and
$P_{j}$ lying on the hyperbolic plane $SNP_{j}$, by the
formula~(\ref{points_distance_coords}) we get:
\begin{equation}\label{|SP_{j}|}
\rho_{j}^2+z_{j}^2+1=2z_{j}\cosinh e_{j}'.
\end{equation}
Now we apply Lemma~\ref{points_distance_lemma} to the vertices $N$
and $P_{j}$:
\begin{equation}\label{|NP_{j}|}
\rho_{j}^2+z_{j}^2+t^2=2tz_{j}\cosinh e_{j}.
\end{equation}
Subtracting~(\ref{|SP_{j}|}) from~(\ref{|NP_{j}|}), under the
assumption that $t\cosinh e_{j}\neq\cosinh e_{j}'$, we get:
\begin{equation}\label{z_{j}}
z_{j}=\frac{t^2-1}{2(t\cosinh e_{j}-\cosinh e_{j}')}.
\end{equation}
Also, taking into account~(\ref{|SP_{j}|}) and (\ref{z_{j}}), we
obtain:
\begin{equation}\label{rho_j^2_2}
\rho_{j}^2=2z_{j}\cosinh e_{j}'-z_{j}^2-1=\frac{(t^2-1)\cosinh
e_{j}'}{(t\cosinh e_{j}-\cosinh e_{j}')}-\frac{(t^2-1)^2}{4(t\cosinh
e_{j}-\cosinh e_{j}')^2}-1.
\end{equation}

Let $\rho_{j,j+1}$ denote the Euclidean distance between the points
$\widetilde{P}_{j}$ and $\widetilde{P}_{j+1}$, $j=1,...,V$ (here and
below $\rho_{V,V+1}\stackrel{\mathrm{def}}{=}\rho_{V,1}$). Applying
Lemma~\ref{points_distance_lemma} to the vertices $P_{j}$ and
$P_{j+1}$, we get:
\begin{equation}\label{rho^2_j_j+1}
\rho_{j,j+1}^2=2z_{j}z_{j+1}\cosinh e_{j,j+1}-z_{j}^2-z_{j+1}^2.
\end{equation}
By the Pythagorean theorem $\rho_{j,j+1}$ is related to the
Cartesian coordinates of $\widetilde{P}_{j}$ and
$\widetilde{P}_{j+1}$ by the formula
\begin{equation}\label{rho_j_j+1_coords}
\rho_{j,j+1}=\sqrt{(x_{j+1}-x_{j})^2+(y_{j+1}-y_{j})^2}.
\end{equation}
By~(\ref{rho_j_theta_j}) the equation~(\ref{rho_j_j+1_coords})
reduces to:
\begin{equation*}\label{rho^2_j_j+1_coords_2}
\rho_{j,j+1}^2=(x_{j}^2+y_{j}^2)+(x_{j+1}^2+y_{j+1}^2)-2(x_{j}x_{j+1}+y_{j}y_{j+1})=
\rho_{j}^2+\rho_{j+1}^2-2(x_{j}x_{j+1}+y_{j}y_{j+1}).
\end{equation*}
Thus, taking into account~(\ref{rho_j^2_2}) and~(\ref{rho^2_j_j+1}),
the expression $x_{j}x_{j+1}+y_{j}y_{j+1}$, which is a part of
$G_{j,j+1}$ from~(\ref{full angle_functions}), is related to the
lengths of edges of $\mathscr{P}$ by the formula
\begin{equation}\label{Re_G_3}
x_{j}x_{j+1}+y_{j}y_{j+1}=\frac{\rho_{j}^2+\rho_{j+1}^2-\rho_{j,j+1}^2}{2}=z_{j}\cosinh
e_{j}'+z_{j+1}\cosinh e_{j+1}'-z_{j}z_{j+1}\cosinh e_{j,j+1}-1.
\end{equation}
Substituting~(\ref{z_{j}}) in~(\ref{Re_G_3}) we get:
\begin{equation*}\label{Re_G_via_edge_lengths_1}
x_{j}x_{j+1}+y_{j}y_{j+1}=\frac{1}{2}\bigg(\frac{(t^2-1)\cosinh
e_{j}'}{(t\cosinh e_{j}-\cosinh e_{j}')}+\frac{(t^2-1)\cosinh
e_{j+1}'}{(t\cosinh e_{j+1}-\cosinh e_{j+1}')}-
\end{equation*}
\begin{equation}\label{Re_G_via_edge_lengths_2}
-\frac{(t^2-1)^2\cosinh e_{j,j+1}}{2(t\cosinh e_{j}-\cosinh
e_{j}')(t\cosinh e_{j+1}-\cosinh e_{j+1}')}-2\bigg).
\end{equation}

Let us now express $x_{j}y_{j+1}-y_{j}x_{j+1}$, which is also a part
of $G_{j,j+1}$, in terms of the length of edges of $\mathscr{P}$.

According to~(\ref{exp(i_theta)}) we know that
\begin{equation}\label{Re_and_Im_G_via_theta}
\cos\theta_{j,j+1}=
\frac{x_{j}x_{j+1}+y_{j}y_{j+1}}{\rho_{j}\rho_{j+1}}\quad\mbox{and}\quad
\sin\theta_{j,j+1}=\frac{x_{j}y_{j+1}-y_{j}x_{j+1}}{\rho_{j}\rho_{j+1}}.
\end{equation}
Note that by definition~(\ref{rho_j_theta_j}), $\rho_{j}>0$,
$j=1,...,V$.

By the Pythagorean trigonometric identity, the formula
\begin{equation}\label{sin_theta_via_cos_theta}
\sin\theta_{j,j+1}=\sigma_{j,j+1}\sqrt{1-\cos^2\theta_{j,j+1}}
\end{equation}
holds true, where $\sigma_{j,j+1}=1$ if $\sin\theta_{j,j+1}\geq0$,
and $\sigma_{j,j+1}=-1$ if $\sin\theta_{j,j+1}<0$ (remind that
$\theta_{j,j+1}$ is determined in~(\ref{theta_j_j+1})).
Then~(\ref{Re_and_Im_G_via_theta})
and~(\ref{sin_theta_via_cos_theta}) imply
\begin{equation*}\label{Im_G_via_theta_1}
x_{j}y_{j+1}-y_{j}x_{j+1}=\rho_{j}\rho_{j+1}\sin\theta_{j,j+1}=
\sigma_{j,j+1}\rho_{j}\rho_{j+1}\sqrt{1-\cos^2\theta_{j,j+1}}=
\end{equation*}
\begin{equation}\label{Im_G_via_theta_2}
=\sigma_{j,j+1}\rho_{j}\rho_{j+1}\sqrt{1-\frac{(x_{j}x_{j+1}+y_{j}y_{j+1})^2}{\rho_{j}^2\rho_{j+1}^2}}=
\sigma_{j,j+1}\sqrt{\rho_{j}^2\rho_{j+1}^2-(x_{j}x_{j+1}+y_{j}y_{j+1})^2}.
\end{equation}
Substituting~(\ref{rho_j^2_2}) and~(\ref{Re_G_via_edge_lengths_2})
in~(\ref{Im_G_via_theta_2}) we get
\begin{equation*}\label{Im_G_via_edge_lengths_1}
x_{j}y_{j+1}-y_{j}x_{j+1}=\sigma_{j,j+1}\bigg[\bigg(\frac{(t^2-1)\cosinh
e_{j}'}{(t\cosinh e_{j}-\cosinh e_{j}')}-\frac{(t^2-1)^2}{4(t\cosinh
e_{j}-\cosinh e_{j}')^2}-1\bigg)\times
\end{equation*}
\begin{equation*}\label{Im_G_via_edge_lengths_2}
\times\bigg(\frac{(t^2-1)\cosinh e_{j+1}'}{(t\cosinh e_{j+1}-\cosinh
e_{j+1}')}-\frac{(t^2-1)^2}{4(t\cosinh e_{j+1}-\cosinh
e_{j+1}')^2}-1\bigg)-\frac{1}{4}\bigg(\frac{(t^2-1)\cosinh
e_{j}'}{(t\cosinh e_{j}-\cosinh e_{j}')}+
\end{equation*}
\begin{equation}\label{Im_G_via_edge_lengths_3}
+\frac{(t^2-1)\cosinh e_{j+1}'}{(t\cosinh e_{j+1}-\cosinh
e_{j+1}')}-\frac{(t^2-1)^2\cosinh e_{j,j+1}}{2(t\cosinh
e_{j}-\cosinh e_{j}')(t\cosinh e_{j+1}-\cosinh
e_{j+1}')}-2\bigg)^2\bigg]^{\frac{1}{2}}.
\end{equation}

Substituting~(\ref{rho_j^2_2}), (\ref{Re_G_via_edge_lengths_2}),
and~(\ref{Im_G_via_edge_lengths_3}) in~(\ref{full angle_coords}) we
obtain the equation of flexibility of a suspension in terms of the
lengths of edges of $\mathscr{P}$.

\section{Proof of the theorem}
\label{Proof_of_the_theorem}

In order to prove the theorem~\ref{sum_of_equator_edges_theorem} we
shall study singular points of the equation of flexibility of a
suspension.

Assume that a nondegenerate suspension $\mathscr{P}$ flexes. Then,
as we have already mentioned in the
section~\ref{The_main_result's_formulation}, the distance $l_{NS}$
between the poles of $\mathscr{P}$ changes during the flex. Let
$t\stackrel{\mathrm{def}}{=}e^{l_{NS}}$ be the parameter of the flex
of $\mathscr{P}$. The identity~(\ref{full angle_functions}) holds
true at every moment $t$ of the flex, as the values of the
expressions $F_{j,j+1}$, $G_{j,j+1}$, $\rho_{j}^2$, $j=1,...,V$,
which make part~(\ref{full angle_functions}), vary as $t$ changes.
Here the functions
$G_{j,j+1}(t)=[x_{j}x_{j+1}+y_{j}y_{j+1}](t)+i[x_{j}y_{j+1}-y_{j}x_{j+1}](t)$
and $\rho_{j}^2(t)$, $j=1,...,V$, are determined
in~(\ref{rho_j^2_2}), ~(\ref{Re_G_via_edge_lengths_2})
and~(\ref{Im_G_via_edge_lengths_3}).

Assume now that for some $j\in\{1,...,V\}$ the dihedral angle
$\theta_{j,j+1}(t)$ remains constant (the value of
$\theta_{j,j+1}(t)$ can also be equal to zero) as $t$ changes. In
this case the length of the edge $N(t)S(t)$ of the tetrahedron
$N(t)S(t)P_{j}(t)P_{j+1}(t)$ must be constant as well (all other
edges of the tetrahedron are also the edges of $\mathscr{P}_t$,
therefore there lengths are fixed), i.e. the value of $t$ does not
change. As we mentioned in the
section~\ref{The_main_result's_formulation}, in this case $\mathcal
P$ can not be flexible. Thus we have the contradiction. Therefore,
the values of the angles $\theta_{j,j+1}(t)$, $j=1,...,V$, change
continuously during the flex. Hence, there exists such an interval
$(t_1,t_2)$ that for all $t\in(t_1,t_2)$ it is true that
$\theta_{j,j+1}(t)\neq0$ for every $j\in\{1,...,V\}$.

We extend both sides of the equation of flexibility~(\ref{full
angle_functions}) as functions in $t$ on the whole complex plane
$\mathbb{C}$. By the theorem on the uniqueness of the analytic
function \cite{Bitsadze1984}, the expression (\ref{full
angle_functions}) remains valid.

Analytic functions $F_{j,j+1}(t)$, $j=1,...,V$, have a finite number
of algebraic singular points. Without loss of generality we can
assume that none of these points lies in the interval $(t_1,t_2)$.
For every $F_{j,j+1}(t)$, $j=1,...,V$, we choose a single-valued
branch $(F_{j,j+1}(t),D)$, where $D\subset\mathbb{C}$ is an
unbounded domain containing $(t_1,t_2)$. Let $\mathcal{W}\subset D$
be a path connecting $t_{0}\in(t_1,t_2)$ and $\infty$, such that
$t_{0}$ is a unique real point of $\mathcal{W}$. Let us calculate
the limit of $F_{j,j+1}(t)$ as $t\rightarrow\infty$ along
$\mathcal{W}$.

Taking into account~(\ref{rho_j^2_2}) we get
\begin{equation}\label{lim_rho_j^2_t^2_2}
\lim\limits_{t\to\infty}\frac{\rho_{j}^2(t)}{t^2}=\lim\limits_{t\to\infty}
\Big[\frac{1}{t^2}\Big(\frac{(t^2-1)\cosinh e_{j}'}{(t\cosinh
e_{j}-\cosinh e_{j}')}-\frac{(t^2-1)^2}{4(t\cosinh e_{j}-\cosinh
e_{j}')^2}-1\Big)\Big]=-\frac{1}{4\cosinh^2 e_{j}}.
\end{equation}
Similarly, from~(\ref{Re_G_via_edge_lengths_2}) we derive that
\begin{equation}\label{lim_Re_G_t^2_3}
\lim\limits_{t\to\infty}\frac{(x_{j}x_{j+1}+y_{j}y_{j+1})(t)}{t^2}=-\frac{\cosinh
e_{j,j+1}}{4\cosinh e_{j}\cosinh e_{j+1}}.
\end{equation}
Also from~(\ref{Im_G_via_theta_2}) and taking into
account~(\ref{lim_rho_j^2_t^2_2}) and~(\ref{lim_Re_G_t^2_3}) we
have:
\begin{equation*}\label{lim_Im_G^2_t^2_1}
\lim\limits_{t\to\infty}\frac{(x_{j}y_{j+1}-y_{j}x_{j+1})^2(t)}{t^4}=
\lim\limits_{t\to\infty}\bigg[\frac{\rho_{j}^2(t)\rho_{j+1}^2(t)-(x_{j}x_{j+1}+y_{j}y_{j+1})^2(t)}{t^4}\bigg]=
\end{equation*}
\begin{equation*}\label{lim_Im_G^2_t^2_2}
=\frac{1}{16\cosinh^2 e_{j}\cosinh^2 e_{j+1}}-\frac{\cosinh^2
e_{j,j+1}}{16\cosinh^2 e_{j}\cosinh^2 e_{j+1}}=\frac{1-\cosinh^2
e_{j,j+1}}{16\cosinh^2 e_{j}\cosinh^2 e_{j+1}}.
\end{equation*}
Hence,
\begin{equation}\label{lim_Im_G_t^2}
\lim\limits_{t\to\infty}\frac{(x_{j}y_{j+1}-y_{j}x_{j+1})(t)}{t^2}=i\sigma_{j,j+1}\frac{\sqrt{\cosinh^2
e_{j,j+1}-1}}{4\cosinh e_{j}\cosinh e_{j+1}},
\end{equation}
where $\sigma_{j,j+1}\in\{+1,-1\}$ is determined by the
single-valued branch $(F_{j,j+1}(t),D)$ and by the path
$\mathcal{W}$.

By definition of $G_{j,j+1}(t)$ and according
to~(\ref{lim_Re_G_t^2_3}) and~(\ref{lim_Im_G_t^2}), we get:
\begin{equation}\label{lim_G_t^2}
\lim\limits_{t\to\infty}\frac{G_{j,j+1}(t)}{t^2}=-\frac{\cosinh
e_{j,j+1}+\sigma_{j,j+1}\sqrt{\cosinh^2 e_{j,j+1}-1}}{4\cosinh
e_{j}\cosinh e_{j+1}}.
\end{equation}
By~(\ref{lim_G_t^2}) and~(\ref{lim_rho_j^2_t^2_2}), the limit of the
left-hand side of~(\ref{full angle_functions}) at
$t\rightarrow\infty$
\begin{equation*}\label{lim_prod_F*}
\lim\limits_{t\to\infty}\prod\limits_{j=1}^{V}F_{j,j+1}(t)=
\lim\limits_{t\to\infty}\prod\limits_{j=1}^{V}\frac{F_{j,j+1}(t)/t^2}{\rho^2_{j}(t)/t^2}=
\prod\limits_{j=1}^{V}\bigg(\cosinh
e_{j,j+1}+\sigma_{j,j+1}\sqrt{{\cosinh}^2 e_{j,j+1}-1}\bigg),
\end{equation*}
and~(\ref{full angle_functions}) at $t\rightarrow\infty$ transforms
to
\begin{equation}\label{flex_eqn_at_infty}
\prod\limits_{j=1}^{V}\bigg(\cosinh
e_{j,j+1}+\sigma_{j,j+1}\sqrt{{\cosinh}^2 e_{j,j+1}-1}\bigg)=1.
\end{equation}

By the following trigonometric identity of hyperbolic geometry,
$\cosinh^2 x-\sinush^2 x=1$, and because $e_{j,j+1}>0$, we have
\begin{equation}\label{sinh_e_j_j+1}
\sqrt{{\cosinh}^2 e_{j,j+1}-1}=\sqrt{{\sinush}^2 e_{j,j+1}}=\sinush
e_{j,j+1}.
\end{equation}
By~(\ref{sinh_e_j_j+1}) the equation~(\ref{flex_eqn_at_infty})
transforms to
\begin{equation}\label{flex_eqn_at_infty_1}
\prod\limits_{j=1}^{V}\big(\cosinh e_{j,j+1}+\sigma_{j,j+1}\sinush
e_{j,j+1}\big)=1.
\end{equation}
By $\cosinh x=\frac{e^{x}+e^{-x}}{2}$ and $\sinush
x=\frac{e^{x}-e^{-x}}{2}$, we have
\begin{equation}\label{left_multiplier_2}
\cosinh e_{j,j+1}+\sigma_{j,j+1}\sinush e_{j,j+1}=\left\{
\begin{array}{rcl}
e^{e_{j,j+1}} & \mbox{ð°-} & \sigma_{j,j+1}=1, \\
e^{-e_{j,j+1}} & \mbox{ð°-} & \sigma_{j,j+1}=-1. \\
\end{array}
\right. =e^{\sigma_{j,j+1}e_{j,j+1}}.
\end{equation}
Substituting~(\ref{left_multiplier_2})
in~(\ref{flex_eqn_at_infty_1}) and taking the logarithm of the
resulting equation, we get~(\ref{flex_eqn_at_infty_3}) $\square$.

The study of the behavior of the equation of flexibility~(\ref{full
angle_functions}) in neighborhoods of other singular points of the
left-hand side of~(\ref{full angle_functions}) did not give us
interesting results: either we were obtaining trivial identities
like $1=1$ (for instance, as $t\rightarrow\pm1$), or the limit of
the left-hand side of the equation of flexibility was too
complicated to distinguish interesting patterns there.

\section{Verification of the necessary flexibility condition of a nondegenerate suspension
for the Bricard-Stachel octahedra in Lobachevsky 3-space}
\label{Verification_of_the_flexibility_condition_on_Bricard's_octahedra}

In 2002 H.~Stachel \cite{Stachel2006} proved the flexibility of the
analogues of the BricardTs octahedra in Lobachevsky 3-space. Let us
verify the validity of the necessary flexibility condition of a
nondegenerate suspension for the Bricard-Stachel octahedra in
Lobachevsky 3-space.

We define an \emph{octahedron} $\mathcal{O}$ as the suspension
$NABCDS$ with the poles $N$ and $S$, and with the vertices of the
equator $A$, $B$, $C$, and $D$. Note that we can consider vertices
$A$ and $C$ as the poles of $\mathcal{O}$ (in this case the
quadrilateral $NDSB$ serves as the equator of $\mathcal{O}$). Also
we can consider vertices $B$ and $D$ as the poles of $\mathcal{O}$
(in this case the quadrilateral $NCSA$ serves as the equator of
$\mathcal{O}$).

\subsection{Bricard-Stachel octahedra of types 1 and 2}
\label{Bricard's_octahedra_of_1st_and_2nd_type}

The procedure of construction of the Bricard-Stachel octahedra of
types 1 and 2 in Lobachevsky 3-space is the same as for the
Bricard's octahedra of types 1 and 2 in Euclidean 3-space
\cite{Stachel2006}, \cite{Alexandrov2010}.

\begin{figure}
\includegraphics[scale=0.5]{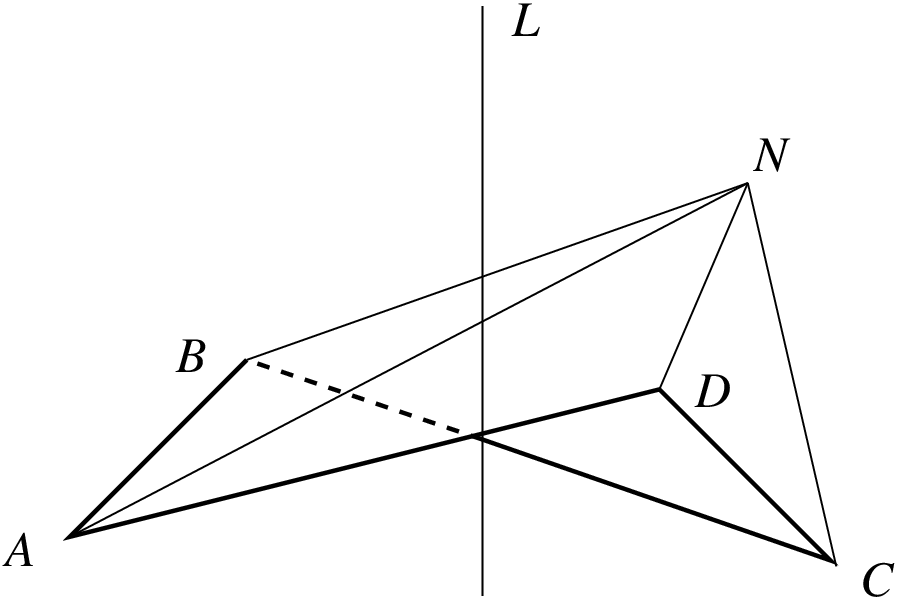} \hfill
\includegraphics[scale=0.5]{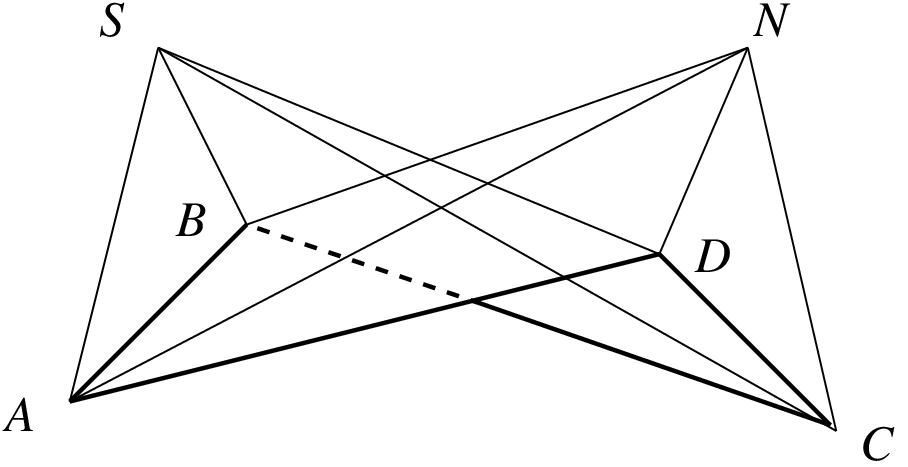} \\
\parbox[t]{7.5cm}{ \caption{The construction of the Bricard-Stachel octahedron of
type 1. Step 1.}\label{bricard1part1}} \hfill
\parbox[t]{7cm}{\caption{The construction of the Bricard-Stachel octahedron of
type 1. Step 2.}\label{bricard1part2}}
\end{figure}

Any \emph{Bricard-Stachel octahedron of type 1} in $\mathbb{H}^3$
can be constructed in the following way. Consider a
disk-homeomorphic piece-wise linear surface $\mathcal{S}$ in
$\mathbb{H}^3$ composed of four triangles $ABN$, $BCN$, $CDN$, and
$DAN$ such that
$\mathrm{d}_{\mathbb{H}^3}(A,B)=\mathrm{d}_{\mathbb{H}^3}(C,D)$ and
$\mathrm{d}_{\mathbb{H}^3}(B,C)=\mathrm{d}_{\mathbb{H}^3}(D,A)$. It
is known that a spatial quadrilateral $ABCD$ which opposite sides
have the same lengths, is symmetric with respect to a line
$\mathcal{L}$ passing through the middle points of its diagonals
$AC$ and $BD$ (see Fig.~\ref{bricard1part1}; for a more precise
analogy with the Euclidean case, in this Figure as well as in the
following Figures we draw polyhedra in the Kleinian model of
Lobachevsky space where lines and planes are intersections of
Euclidean lines and planes with a fixed unit ball). Glue together
$\mathcal{S}$ and its symmetric image with respect to $L$ along
$ABCD$. Denote by $S$ the symmetric image of $N$ under the symmetry
with respect to $L$ (see Fig.~\ref{bricard1part2}). The resulting
polyhedral surface $NABCDS$ with self-intersections is flexible
(because $\mathcal{S}$ is flexible) and combinatorially it is an
octahedron (according to the definition given above). We will call
it a Bricard-Stachel octahedron of type 1. By construction it
follows that
$\mathrm{d}_{\mathbb{H}^3}(A,N)=\mathrm{d}_{\mathbb{H}^3}(C,S)$,
$\mathrm{d}_{\mathbb{H}^3}(B,N)=\mathrm{d}_{\mathbb{H}^3}(D,S)$,
$\mathrm{d}_{\mathbb{H}^3}(C,N)=\mathrm{d}_{\mathbb{H}^3}(A,S)$, and
$\mathrm{d}_{\mathbb{H}^3}(D,N)=\mathrm{d}_{\mathbb{H}^3}(B,S)$.

\begin{figure}
\includegraphics[scale=0.5]{bricard2part1.eps}
\hfill
\includegraphics[scale=0.5]{bricard2part2.eps} \\
\parbox[t]{7.5cm}{ \caption{The construction of the Bricard-Stachel octahedron of
type 2. Step 1.}\label{bricard2part1}} \hfill
\parbox[t]{7cm}{\caption{The construction of the Bricard-Stachel octahedron of
type 2. Step 2.}\label{bricard2part2}}
\end{figure}

Any \emph{Bricard-Stachel octahedron of type 2} in $\mathbb{H}^3$
can be constructed as follows. Consider a disk-homeomorphic
piece-wise linear surface $\mathcal{S}$ in $\mathbb{H}^3$ composed
of four triangles $ABN$, $BCN$, $CDN$, and $DAN$ such that
$\mathrm{d}_{\mathbb{H}^3}(A,B)=\mathrm{d}_{\mathbb{H}^3}(B,C)$ and
$\mathrm{d}_{\mathbb{H}^3}(C,D)=\mathrm{d}_{\mathbb{H}^3}(D,A)$. It
is known that a spatial quadrilateral $ABCD$ which neighbor sides at
the vertices $B$ and $D$ have the same lengths, is symmetric with
respect to a plane $H$ which dissects the dihedral angle between the
half-planes $ABD$ and $CBD$ (see Fig.~\ref{bricard2part1}). Glue
together $\mathcal{S}$ and its symmetric image with respect to $H$
along $ABCD$. Denote by $S$ the symmetric image of $N$ under the
symmetry with respect to $H$ (see Fig.~\ref{bricard1part2}). The
resulting polyhedral surface $NABCDS$ with self-intersections is
flexible (because $\mathcal{S}$ is flexible) and combinatorially it
is an octahedron. We will call it a Bricard-Stachel octahedron of
type 2. By construction it follows that
$\mathrm{d}_{\mathbb{H}^3}(A,N)=\mathrm{d}_{\mathbb{H}^3}(C,S)$,
$\mathrm{d}_{\mathbb{H}^3}(C,N)=\mathrm{d}_{\mathbb{H}^3}(A,S)$,
$\mathrm{d}_{\mathbb{H}^3}(B,N)=\mathrm{d}_{\mathbb{H}^3}(B,S)$, and
$\mathrm{d}_{\mathbb{H}^3}(D,N)=\mathrm{d}_{\mathbb{H}^3}(D,S)$.

It remains to note that for every considered octahedron each of
three its equators has two pairs of edges of the same lengths.
Hence, the theorem~\ref{sum_of_equator_edges_theorem} is valid for
the Bricard-Stachel octahedra of types 1 and 2.

\subsection{Bricard-Stachel octahedra of type 3}
\label{Bricard's_octahedra_of_3rd_type}

There are three subtypes of the Bricard-Stachel octahedra of type 3
in Lobachevsky space \cite{Stachel2006} which construction is based
on circles, horocycles or hypercircles correspondingly. The
procedure of construction is common for all subtypes of the
Bricard-Stachel octahedra of type 3 and it is the same as for the
Bricard's octahedra of type 3 in Euclidean space.

\begin{figure}
\begin{center}
\includegraphics[scale=0.5]{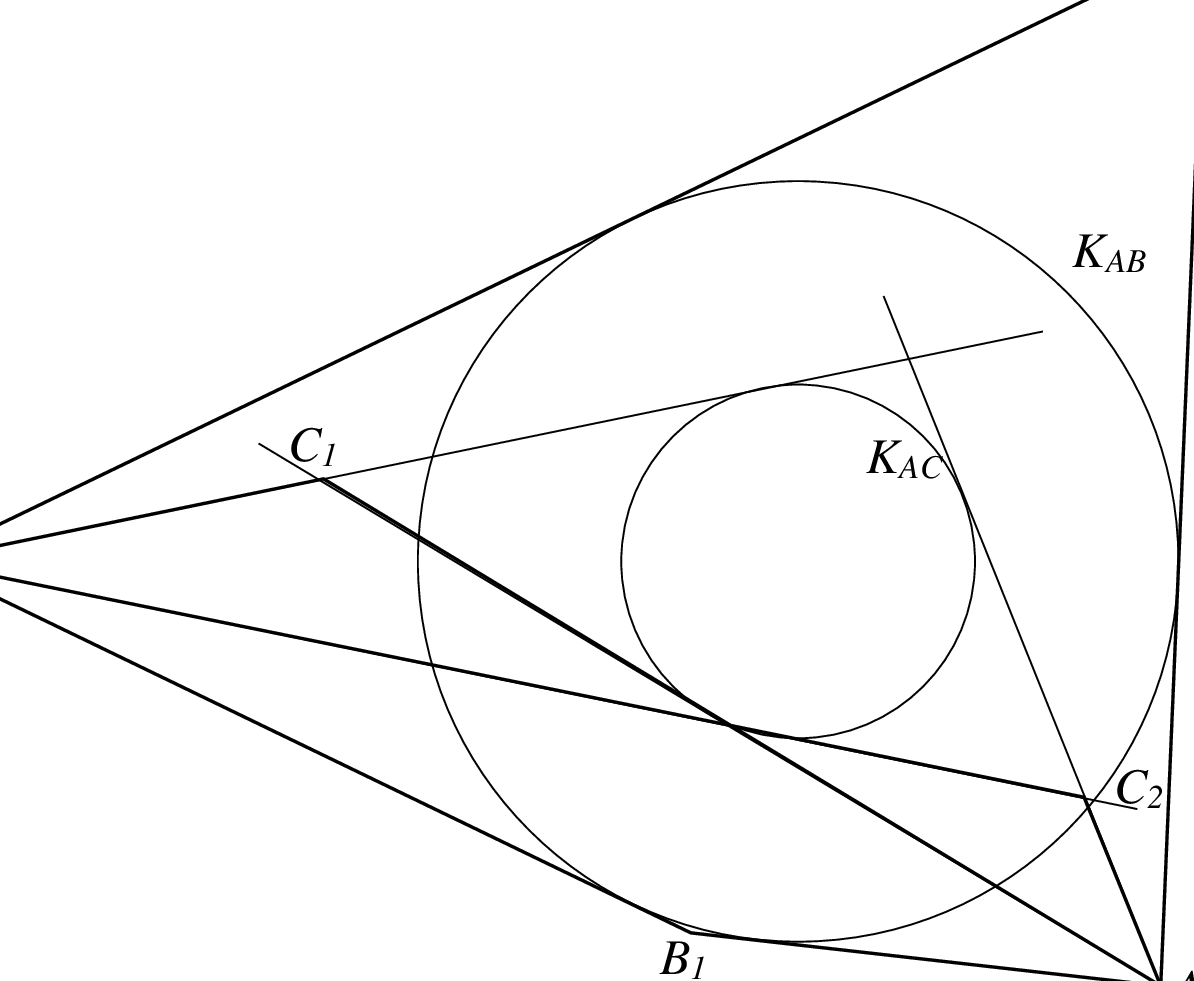}
\caption{The construction of the Bricard-Stachel octahedron of type
3 based on circles. Step 1.}\label{bricard3part1}
\end{center}
\end{figure}

\begin{figure}
\begin{center}
\includegraphics[scale=0.5]{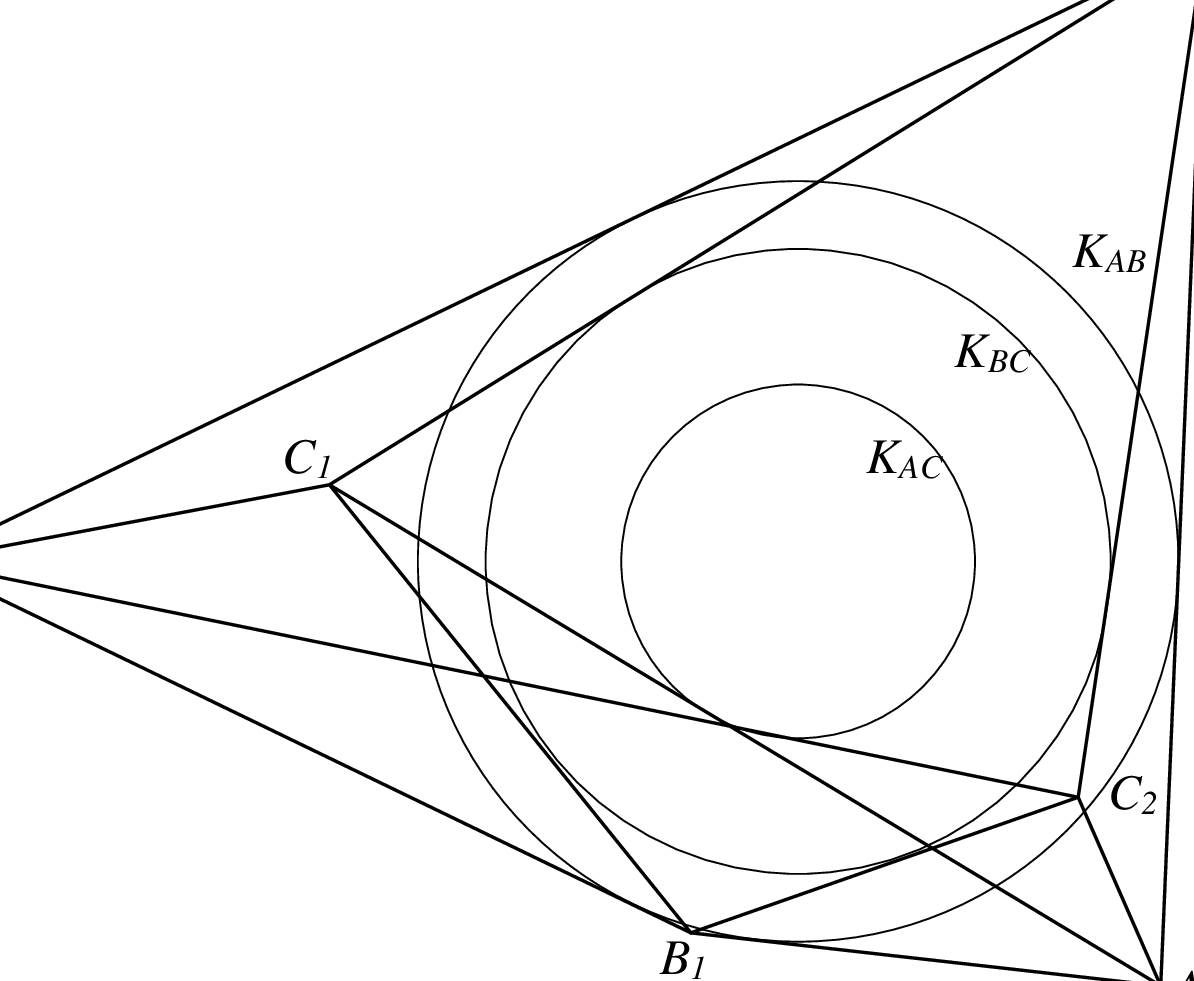}
\caption{The construction of the Bricard-Stachel octahedron of type
3 based on circles. Step 2.}\label{bricard3part2}
\end{center}
\end{figure}

Any \emph{Bricard-Stachel octahedron of type 3} in $\mathbb{H}^3$
can be constructed in the following way. Let $K_{AC}$ and $K_{AB}$
be two different circles (horocycles, hypercircles) in
$\mathbb{H}^2$ with the common center $M$ and let $A_{1}$, $A_{2}$
be two different finite points outside $K_{AC}$ and $K_{AB}$. In
addition, suppose that $K_{AC}$, $K_{AB}$, $A_{1}$ and $A_{2}$ are
taken in such a way that the straight lines tangent to $K_{AB}$ and
passing through $A_{1}$ and $A_{2}$ intersect pairwise in finite
points of $\mathbb{H}^2$ and form a quadrilateral
$A_{1}B_{1}A_{2}B_{2}$ tangent to $K_{AB}$; moreover, that the
straight lines tangent to $K_{AC}$ and passing through $A_{1}$ and
$A_{2}$ intersect pairwise in finite points of $\mathbb{H}^2$ and
form a quadrilateral $A_{1}C_{1}A_{2}C_{2}$ tangent to $K_{AC}$ (see
Fig.~\ref{bricard3part1}; for clarity, we placed circles $K_{AB}$
and $K_{AC}$ so that their common center coincides with the center
of the Kleinian model of Lobachevsky space. In this case $K_{AB}$
and $K_{AC}$ are Euclidean circles as well). A polyhedron
$\mathcal{O}$ with the vertices $A_{i}$, $B_{j}$, $C_{k}$, with the
edges $A_{i}B_{j}$, $A_{i}C_{k}$, $B_{j}C_{k}$, and with the faces
$\triangle A_{i}B_{j}C_{k}$, $i,j,k\in\{1,2\}$, is an octahedron in
the sense of the definition given above (see
Fig.~\ref{bricard3part2}). The following pairs of vertices can serve
as the poles of $\mathcal{O}$: $(A_{1},A_{2})$ with the
corresponding equator $B_{1}C_{1}B_{2}C_{2}$, $(B_{1},B_{2})$ with
the equator $A_{1}C_{1}A_{2}C_{2}$, and $(C_{1},C_{2})$ with the
equator $A_{1}B_{1}A_{2}B_{2}$. Suppose in addition that
$\mathcal{O}$ does not have symmetries. We will call such octahedron
$\mathcal{O}$ a Bricard-Stachel octahedron of type 3.

According to H.~Stachel \cite{Stachel2006}, $\mathcal{O}$ flexes
continuously in $\mathbb{H}^3$. Moreover, $\mathcal{O}$ admits two
flat positions during the flex (we constructed $\mathcal{O}$ in one
of its flat positions). Hence, for every equator of $\mathcal{O}$,
$A_{1}B_{1}A_{2}B_{2}$, $B_{1}C_{1}B_{2}C_{2}$, and
$A_{1}C_{1}A_{2}C_{2}$, all straight lines containing a side of the
equator are tangent to some circle (horocycle, hypercircle) at least
in one flat position of $\mathcal{O}$. Using this fact, we will
prove that the theorem~\ref{sum_of_equator_edges_theorem} is valid
for the Bricard-Stachel octahedra of type 3. We have to consider
three possible cases: when an equator of $\mathcal{O}$ is tangent to
a circle, to a horocycle, or to a hypercircle in $\mathbb{H}^2$.
Here we study the most common situation when any three vertices of
an equator of a flexible octahedron in its flat position do not lie
on a straight line.

\subsubsection{An equator of a Bricard-Stachel octahedron of type 3
is tangent to a circle in $\mathbb{H}^2$}
\label{Bricard's_octahedra_of_3rd_type_tangent to_a_circle}

Let $M$ be the center of the circle $K_{AB}$ with the radius $R$ in
$\mathbb{H}^2$ and let all straight lines containing a side of the
quadrilateral $A_{1}B_{1}A_{2}B_{2}$ are tangent to $K_{AB}$. Let us
draw the segments $MP_{1}$, $MP_{2}$, $MP_{3}$, $MP_{4}$ connecting
$M$ with the straight lines $A_{1}B_{2}$, $A_{2}B_{2}$,
$A_{2}B_{1}$, $A_{1}B_{1}$ and perpendicular to the corresponding
lines. By construction,
$\mathrm{d}_{\mathbb{H}^2}(M,P_{1})=\mathrm{d}_{\mathbb{H}^2}(M,P_{2})=
\mathrm{d}_{\mathbb{H}^2}(M,P_{3})=\mathrm{d}_{\mathbb{H}^2}(M,P_{4})=R$.

By the Pythagorean theorem for Lobachevsky space \cite{Vinberg1988}
applied to $\triangle A_{1}MP_{1}$ and $\triangle A_{1}MP_{4}$, we
obtain: $\cosinh \mathrm{d}_{\mathbb{H}^2}(A_{1},P_{1})=\cosinh
\mathrm{d}_{\mathbb{H}^2}(A_{1},P_{4})=\cosinh
\mathrm{d}_{\mathbb{H}^2}(A_{1},M)/\cosinh R$. Then
$a\stackrel{\mathrm{def}}{=}\mathrm{d}_{\mathbb{H}^2}(A_{1},P_{1})=\mathrm{d}_{\mathbb{H}^2}(A_{1},P_{4})$.
Similarly we get:
$b\stackrel{\mathrm{def}}{=}\mathrm{d}_{\mathbb{H}^2}(B_{2},P_{1})=\mathrm{d}_{\mathbb{H}^2}(B_{2},P_{2})$,
$c\stackrel{\mathrm{def}}{=}\mathrm{d}_{\mathbb{H}^2}(A_{2},P_{2})=\mathrm{d}_{\mathbb{H}^2}(A_{2},P_{3})$,
and
$d\stackrel{\mathrm{def}}{=}\mathrm{d}_{\mathbb{H}^2}(B_{1},P_{3})=\mathrm{d}_{\mathbb{H}^2}(B_{1},P_{4})$.

If the circle $K_{AB}$ is inscribed in the quadrilateral
$A_{1}B_{1}A_{2}B_{2}$ (see Fig.~\ref{bricard3part1}), then
\linebreak $\mathrm{d}_{\mathbb{H}^2}(A_{1},B_{2})=a+b$,
$\mathrm{d}_{\mathbb{H}^2}(A_{2},B_{2})=b+c$,
$\mathrm{d}_{\mathbb{H}^2}(A_{2},B_{1})=c+d$,
$\mathrm{d}_{\mathbb{H}^2}(A_{1},B_{1})=a+d$, and the identity
\begin{equation}\label{equator_circle_+-+-}
\mathrm{d}_{\mathbb{H}^2}(A_{1},B_{2})-\mathrm{d}_{\mathbb{H}^2}(A_{2},B_{2})+
\mathrm{d}_{\mathbb{H}^2}(A_{1},B_{1})-\mathrm{d}_{\mathbb{H}^2}(A_{1},B_{1})=0
\end{equation}
holds true.

If the circle $K_{AB}$ is tangent to the quadrilateral
$A_{1}B_{1}A_{2}B_{2}$ externally (this case corresponds to the
quadrilateral $A_{1}C_{1}A_{2}C_{2}$ and to the circle $K_{AC}$ in
the Fig.~\ref{bricard3part1}), then
$\mathrm{d}_{\mathbb{H}^2}(A_{1},B_{2})=a-b$,
$\mathrm{d}_{\mathbb{H}^2}(A_{2},B_{2})=b+c$,
$\mathrm{d}_{\mathbb{H}^2}(A_{2},B_{1})=c-d$,
$\mathrm{d}_{\mathbb{H}^2}(A_{1},B_{1})=a+d$, and the identity
\begin{equation}\label{equator_circle_++--}
\mathrm{d}_{\mathbb{H}^2}(A_{1},B_{2})+\mathrm{d}_{\mathbb{H}^2}(A_{2},B_{2})-
\mathrm{d}_{\mathbb{H}^2}(A_{1},B_{1})-\mathrm{d}_{\mathbb{H}^2}(A_{1},B_{1})=0
\end{equation}
holds true.

By~(\ref{equator_circle_+-+-}) and~(\ref{equator_circle_++--}), the
theorem~\ref{sum_of_equator_edges_theorem} is valid for any equator
of a Bricard-Stachel octahedron of type 3 tangent to a circle in at
least one of its flat positions.

\subsubsection{An equator of a Bricard-Stachel octahedron of type 3
is tangent to a horocycle in $\mathbb{H}^2$}
\label{Bricard's_octahedra_of_3rd_type_tangent to_a_horocycle}

Let us consider the Poincar\'e upper half-plane model of Lobachevsky
plane $\mathbb{H}^2$ with the coordinates $(\rho,z)$ (i.e., with the
metric given by the formula $ds^2=\frac{d\rho^2+dz^2}{z^2}$).
Without loss of generality we can assume that the center of the
horocycle tangent to the equator of a Bricard-Stachel octahedron
$\mathcal{O}$ of type 3, coincides with the (unique) point $\infty$
at infinity of $\mathbb{H}^2$ which does not lie on the Euclidean
line $z=0$. We denote the family of such horocycles by
$K=\{\rho=R|R>0\}$. Let $K_{R}\in K$ and let $A_1=({\rho}_{A_1},
z_{A_1})$ and $A_2=({\rho}_{A_2}, z_{A_2})$ be two opposite vertices
of $\mathcal{O}$, such that the straight line (in $\mathbb{H}^2$)
passing through $A_1$ and $A_2$ is not tangent to $K_R$. All the
vertices of $\mathcal{O}$ are located outside $K_R$, hence
$z_{A_1}<R$ and $z_{A_2}<R$. We will construct all possible
quadrangles tangent to $K_R$ with the opposite vertices $A_1$ and
$A_2$, i.e., all quadrangles that can serve as equators of
$\mathcal{O}$. Then we will verify the validity of the
theorem~\ref{sum_of_equator_edges_theorem} for such quadrangles.

Let $T=({\rho}_{T}, z_{T})$ be a point in $\mathbb{H}^2$ and let
$\Lambda$ be a straight line in $\mathbb{H}^2$ passing through $T$
which is realized in the Poincar\'e upper half-plane as the
Euclidean demi-circle with the radius
$\sqrt{({\rho}_{T}-{\rho}_{T,\Lambda})^2+z_{T}^2}$ and with the
center $O^{T}_{\Lambda}=({\rho}_{T,\Lambda},0)$. Then the angle
$\varphi_{T}^{\Lambda}\stackrel{\mathrm{def}}{=}\angle
TO^{T}_{\Lambda}\rho\in(0,\pi)$ determines uniquely a position of
$T$ on $\Lambda$.

\begin{remark}\label{r1} For every finite point $T=({\rho}_{T},
z_{T})$, $z_{T}<R$, there exist precisely two straight lines
${\Lambda}_{l}^{T}$ and ${\Lambda}_{r}^{T}$ tangent to the horocycle
$K_R$ and containing $T$. They are realized in the Poincar\'e upper
half-plane as the Euclidean demi-circles with the radius $R$ and
with the centers $O^{T}_{l}=({\rho}_{T,l},0)$ and
$O^{T}_{r}=({\rho}_{T,r},0)$,
${\rho}_{T,l}\leq{\rho}_{T}\leq{\rho}_{T,r}$. The angles
$\varphi_{T}^{l}\stackrel{\mathrm{def}}{=}\angle TO^{T}_{l}\rho$ and
$\varphi_{T}^{r}\stackrel{\mathrm{def}}{=}\angle TO^{T}_{r}\rho$
serve as the coordinates of $T$ on ${\Lambda}_{l}^{T}$ and
${\Lambda}_{r}^{T}$ correspondingly. Then, by construction, we get:
$\varphi_{T}^{r}=\pi-\varphi_{T}^{l}$. Hence,
\begin{equation}\label{cos(right)=-cos(left)}
\cos\varphi_{T}^{r}=-\cos\varphi_{T}^{l}.
\end{equation}
\end{remark}

According to the remark~\ref{r1}, there are two straight lines,
${\Lambda}_{l}^{A_1}$, and ${\Lambda}_{r}^{A_1}$, passing through
$A_1$ and tangent to $K_R$, which are realised in $\mathbb{H}^2$ as
the Euclidean demi-circles with the radius $R$ and with the centers
$O^{A_1}_{l}=({\rho}_{A_1,l},0)$, $O^{A_1}_{r}=({\rho}_{A_1,r},0)$,
${\rho}_{A_1,l}\leq{\rho}_{A_1}\leq{\rho}_{A_1,r}$. The angles
$\varphi_{A_1}^{{\Lambda}_{l}^{A_1}}\stackrel{\mathrm{def}}{=}\angle
A_1 O^{A_1}_{l}\rho$,
$\varphi_{A_1}^{{\Lambda}_{r}^{A_1}}\stackrel{\mathrm{def}}{=}\angle
A_1 O^{A_1}_{r} \rho$ serve as the coordinates of $A_1$ on
${\Lambda}_{l}^{A_1}$ and ${\Lambda}_{r}^{A_1}$ correspondingly.
Moreover,
\begin{equation}\label{phi_A_1}
\cos\varphi_{A_1}^{{\Lambda}_{r}^{A_1}}=-\cos\varphi_{A_1}^{{\Lambda}_{l}^{A_1}}.
\end{equation}

Similarly, there are two straight lines, ${\Lambda}_{l}^{A_2}$, and
${\Lambda}_{r}^{A_2}$, passing through $A_2$ and tangent to $K_R$,
which are realised in $\mathbb{H}^2$ as the Euclidean demi-circles
with the radius $R$ and with the centers
$O^{A_2}_{l}=({\rho}_{A_2,l},0)$, $O^{A_2}_{r}=({\rho}_{A_2,r},0)$,
${\rho}_{A_2,l}\leq{\rho}_{A_2}\leq{\rho}_{A_2,r}$. The angles
$\varphi_{A_2}^{{\Lambda}_{l}^{A_2}}\stackrel{\mathrm{def}}{=}\angle
A_2 O^{A_2}_{l}\rho$,
$\varphi_{A_2}^{{\Lambda}_{r}^{A_2}}\stackrel{\mathrm{def}}{=}\angle
A_2 O^{A_2}_{r} \rho$ serve as the coordinates of $A_2$ on
${\Lambda}_{l}^{A_2}$ and ${\Lambda}_{r}^{A_2}$ correspondingly.
Moreover,
\begin{equation}\label{phi_A_2}
\cos\varphi_{A_2}^{{\Lambda}_{r}^{A_2}}=-\cos\varphi_{A_2}^{{\Lambda}_{l}^{A_2}}.
\end{equation}

Suppose that ${\Lambda}_{l}^{A_1}$ and ${\Lambda}_{l}^{A_2}$
intersect at a point $B_1$. Then the angles
$\varphi_{B_1}^{{\Lambda}_{l}^{A_1}}\stackrel{\mathrm{def}}{=}\angle
B_1 O^{A_1}_{l}\rho$,
$\varphi_{B_1}^{{\Lambda}_{l}^{A_2}}\stackrel{\mathrm{def}}{=}\angle
B_1 O^{A_2}_{l} \rho$ serve as the coordinates of $B_1$ on
${\Lambda}_{l}^{A_1}$ and ${\Lambda}_{l}^{A_2}$ correspondingly.
Moreover,
\begin{equation}\label{phi_B_1}
\cos\varphi_{B_1}^{{\Lambda}_{l}^{A_2}}=-\cos\varphi_{B_1}^{{\Lambda}_{l}^{A_1}}.
\end{equation}

Also suppose that ${\Lambda}_{r}^{A_1}$ and ${\Lambda}_{r}^{A_2}$
intersect at a point $B_2$. Then the angles
$\varphi_{B_2}^{{\Lambda}_{r}^{A_1}}\stackrel{\mathrm{def}}{=}\angle
B_2 O^{A_1}_{r}\rho$,
$\varphi_{B_2}^{{\Lambda}_{r}^{A_2}}\stackrel{\mathrm{def}}{=}\angle
B_2 O^{A_2}_{r} \rho$ serve as the coordinates of $B_2$ on
${\Lambda}_{r}^{A_1}$ and ${\Lambda}_{r}^{A_2}$ correspondingly.
Moreover,
\begin{equation}\label{phi_B_2}
\cos\varphi_{B_2}^{{\Lambda}_{r}^{A_2}}=-\cos\varphi_{B_2}^{{\Lambda}_{r}^{A_1}}.
\end{equation}

Let the straight lines ${\Lambda}_{r}^{A_1}$ and
${\Lambda}_{l}^{A_2}$ intersect at a point $C_1$. Then the angles
$\varphi_{C_1}^{{\Lambda}_{r}^{A_1}}\stackrel{\mathrm{def}}{=}\angle
C_1 O^{A_1}_{r}\rho$,
$\varphi_{C_1}^{{\Lambda}_{l}^{A_2}}\stackrel{\mathrm{def}}{=}\angle
C_1 O^{A_2}_{l} \rho$ serve as the coordinates of $C_1$ on
${\Lambda}_{r}^{A_1}$ and ${\Lambda}_{l}^{A_2}$ correspondingly.
Moreover,
\begin{equation}\label{phi_C_1}
\cos\varphi_{C_1}^{{\Lambda}_{l}^{A_2}}=-\cos\varphi_{C_1}^{{\Lambda}_{r}^{A_1}}.
\end{equation}

Also, let the straight lines ${\Lambda}_{l}^{A_1}$ and
${\Lambda}_{r}^{A_2}$ intersect at a point $C_2$. Then the angles
$\varphi_{C_2}^{{\Lambda}_{l}^{A_1}}\stackrel{\mathrm{def}}{=}\angle
C_2 O^{A_1}_{l}\rho$,
$\varphi_{C_2}^{{\Lambda}_{r}^{A_2}}\stackrel{\mathrm{def}}{=}\angle
C_2 O^{A_2}_{r} \rho$ serve as the coordinates of $C_2$ on
${\Lambda}_{l}^{A_1}$ and ${\Lambda}_{r}^{A_2}$ correspondingly.
Moreover,
\begin{equation}\label{phi_C_2}
\cos\varphi_{C_2}^{{\Lambda}_{r}^{A_2}}=-\cos\varphi_{C_2}^{{\Lambda}_{l}^{A_1}}.
\end{equation}

By construction, the quadrangles $A_1 B_1 A_2 B_2$ and $A_1 C_1 A_2
C_2$ are tangent to $K_R$, and the points $A_1$, $A_2$ are opposite
vertices of each of these quadrangles. In order to verify the
validity of the theorem~\ref{sum_of_equator_edges_theorem} for the
flexible octahedra with the equator $A_1 B_1 A_2 B_2$ or $A_1 C_1
A_2 C_2$ we need to prove the following easy statement.

\begin{lemma}\label{points_distance_via_angles_lemma}
Given a Poincar\'e upper half-plane $\mathbb{H}^2$ with the
coordinates $(\rho,z)$ \emph{(}i.e., with the metric given by the
formula $ds^2=\frac{d\rho^2+dz^2}{z^2}$\emph{)}. Let $A$ and $B$ be
points on the straight line $\Lambda$ realized in $\mathbb{H}^2$ as
the Euclidean demi-circle with the raduis $R$ and with the center
$O_{\Lambda}=({\rho}_{O_{\Lambda}},0)$, and let the angles
$\varphi_A\stackrel{\mathrm{def}}{=}\angle AO_{\Lambda}\rho$,
$\varphi_B\stackrel{\mathrm{def}}{=}\angle BO_{\Lambda}\rho$ serve
as the coordinates of $A$ and $B$ correspondingly on $\Lambda$. Also
we assume that $0<\varphi_A\leq\phi_B<\pi$. Then the distance
between $A$ and $B$ is calculated as follows:
\begin{equation}\label{points_distance_coords_via_angles}
\mathrm{d}_{\mathbb{H}^2}(A,B)=\frac{1}{2}\ln\Big[\Big(\frac{1+\cos\varphi_{A}}{1+\cos\varphi_{B}}\Big)
\Big(\frac{1-\cos\varphi_{B}}{1-\cos\varphi_{A}}\Big)\Big].
\end{equation}
\end{lemma}

\textbf{Proof.} The hyperbolic segment $\Lambda_{AB}$ connecting the
points $A$ and $B$ is specified parametrically by the formulas
$\Lambda_{AB}(t):(\rho(\varphi),z(\varphi))$,
$\varphi\in[\varphi_{A},\varphi_{B}]$, where
$\rho(\varphi)={\rho}_{O_{\Lambda}}+R\cos\varphi$,
$z(\varphi)=R\sin\varphi$, $A=\Lambda_{AB}(\varphi_{A})$,
$B=\Lambda_{AB}(\varphi_{B})$. The direct calculation shows that the
lengths of $\Lambda_{AB}$ is equal to the right-hand side
of~(\ref{points_distance_coords_via_angles}). $\square$

By Lemma~\ref{points_distance_via_angles_lemma}, the lengths of the
edges of the quadrilateral $A_1 B_1 A_2 B_2$ are calculated as
follows:
\begin{equation}\label{length_of_A_1_B_1_horocycle}
\mathrm{d}_{\mathbb{H}^2}(A_1,B_1)=\frac{1}{2}\ln\Bigg[\Bigg(\frac{1+\cos\varphi_{A_1}^{{\Lambda}_{l}^{A_1}}}
{1+\cos\varphi_{B_1}^{{\Lambda}_{l}^{A_1}}}\Bigg)
\Bigg(\frac{1-\cos\varphi_{B_1}^{{\Lambda}_{l}^{A_1}}}{1-\cos\varphi_{A_1}^{{\Lambda}_{l}^{A_1}}}\Bigg)\Bigg],
\end{equation}
\begin{equation}\label{length_of_A_2_B_1_horocycle}
\mathrm{d}_{\mathbb{H}^2}(A_2,B_1)=\frac{1}{2}\ln\Bigg[\Bigg(\frac{1+\cos\varphi_{A_2}^{{\Lambda}_{l}^{A_2}}}
{1+\cos\varphi_{B_1}^{{\Lambda}_{l}^{A_2}}}\Bigg)
\Bigg(\frac{1-\cos\varphi_{B_1}^{{\Lambda}_{l}^{A_2}}}{1-\cos\varphi_{A_2}^{{\Lambda}_{l}^{A_2}}}\Bigg)\Bigg],
\end{equation}
\begin{equation}\label{length_of_B_2_A_1_horocycle}
\mathrm{d}_{\mathbb{H}^2}(B_2,A_1)=\frac{1}{2}\ln\Bigg[\Bigg(\frac{1+\cos\varphi_{B_2}^{{\Lambda}_{r}^{A_1}}}
{1+\cos\varphi_{A_1}^{{\Lambda}_{r}^{A_1}}}\Bigg)
\Bigg(\frac{1-\cos\varphi_{A_1}^{{\Lambda}_{r}^{A_1}}}{1-\cos\varphi_{B_2}^{{\Lambda}_{r}^{A_1}}}\Bigg)\Bigg],
\end{equation}
\begin{equation}\label{length_of_B_2_A_2_horocycle}
\mathrm{d}_{\mathbb{H}^2}(B_2,A_2)=\frac{1}{2}\ln\Bigg[\Bigg(\frac{1+\cos\varphi_{B_2}^{{\Lambda}_{r}^{A_2}}}
{1+\cos\varphi_{A_2}^{{\Lambda}_{r}^{A_2}}}\Bigg)
\Bigg(\frac{1-\cos\varphi_{A_2}^{{\Lambda}_{r}^{A_2}}}{1-\cos\varphi_{B_2}^{{\Lambda}_{r}^{A_2}}}\Bigg)\Bigg].
\end{equation}
Then, by (\ref{phi_A_1})---(\ref{phi_B_2}), we get:
\begin{equation}\label{equator_A_1_B_1_A_2_B_2_horocycle}
\mathrm{d}_{\mathbb{H}^2}(A_1,B_1)+\mathrm{d}_{\mathbb{H}^2}(A_2,B_1)-
\mathrm{d}_{\mathbb{H}^2}(B_2,A_1)-\mathrm{d}_{\mathbb{H}^2}(B_2,A_2)=0.
\end{equation}

By Lemma~\ref{points_distance_via_angles_lemma}, the lengths of the
edges of the quadrilateral $A_1 C_1 A_2 C_2$ are calculated as
follows:
\begin{equation}\label{length_of_C_1_A_1_horocycle}
\mathrm{d}_{\mathbb{H}^2}(C_1,A_1)=\frac{1}{2}\ln\Bigg[\Bigg(\frac{1+\cos\varphi_{C_1}^{{\Lambda}_{r}^{A_1}}}
{1+\cos\varphi_{A_1}^{{\Lambda}_{r}^{A_1}}}\Bigg)
\Bigg(\frac{1-\cos\varphi_{A_1}^{{\Lambda}_{r}^{A_1}}}{1-\cos\varphi_{C_1}^{{\Lambda}_{r}^{A_1}}}\Bigg)\Bigg],
\end{equation}
\begin{equation}\label{length_of_C_2_A_1_horocycle}
\mathrm{d}_{\mathbb{H}^2}(C_2,A_1)=\frac{1}{2}\ln\Bigg[\Bigg(\frac{1+\cos\varphi_{C_2}^{{\Lambda}_{l}^{A_1}}}
{1+\cos\varphi_{A_1}^{{\Lambda}_{l}^{A_1}}}\Bigg)
\Bigg(\frac{1-\cos\varphi_{A_1}^{{\Lambda}_{l}^{A_1}}}{1-\cos\varphi_{C_2}^{{\Lambda}_{l}^{A_1}}}\Bigg)\Bigg],
\end{equation}
\begin{equation}\label{length_of_A_2_C_1_horocycle}
\mathrm{d}_{\mathbb{H}^2}(A_2,C_1)=\frac{1}{2}\ln\Bigg[\Bigg(\frac{1+\cos\varphi_{A_2}^{{\Lambda}_{l}^{A_2}}}
{1+\cos\varphi_{C_1}^{{\Lambda}_{l}^{A_2}}}\Bigg)
\Bigg(\frac{1-\cos\varphi_{C_1}^{{\Lambda}_{l}^{A_2}}}{1-\cos\varphi_{A_2}^{{\Lambda}_{l}^{A_2}}}\Bigg)\Bigg],
\end{equation}
\begin{equation}\label{length_of_A_2_C_2_horocycle}
\mathrm{d}_{\mathbb{H}^2}(A_2,C_2)=\frac{1}{2}\ln\Bigg[\Bigg(\frac{1+\cos\varphi_{A_2}^{{\Lambda}_{r}^{A_2}}}
{1+\cos\varphi_{C_2}^{{\Lambda}_{r}^{A_2}}}\Bigg)
\Bigg(\frac{1-\cos\varphi_{C_2}^{{\Lambda}_{r}^{A_2}}}{1-\cos\varphi_{A_2}^{{\Lambda}_{r}^{A_2}}}\Bigg)\Bigg].
\end{equation}
By (\ref{phi_A_1}), (\ref{phi_A_2}), (\ref{phi_C_1}), and
(\ref{phi_C_2}), it is easy to verify that
\begin{equation}\label{equator_A_1_C_1_A_2_C_2_horocycle}
\mathrm{d}_{\mathbb{H}^2}(C_2,A_1)+\mathrm{d}_{\mathbb{H}^2}(C_1,A_1)-
\mathrm{d}_{\mathbb{H}^2}(A_2,C_1)-\mathrm{d}_{\mathbb{H}^2}(A_2,C_2)=0.
\end{equation}

According to~(\ref{equator_A_1_B_1_A_2_B_2_horocycle})
and~(\ref{equator_A_1_C_1_A_2_C_2_horocycle}), the
theorem~\ref{sum_of_equator_edges_theorem} is valid for any equator
of a Bricard-Stachel octahedron of type 3 tangent to a horocycle in
at least one of its flat positions.

\subsubsection{An equator of a Bricard-Stachel octahedron of type 3
is tangent to a hypercircle in $\mathbb{H}^2$}
\label{Bricard's_octahedra_of_3rd_type_tangent to_a_hypercircle}

Let us consider the Poincar\'e upper half-plane model of Lobachevsky
plane $\mathbb{H}^2$ with the coordinates $(\rho,z)$ (i.e., with the
metric given by the formula $ds^2=\frac{d\rho^2+dz^2}{z^2}$).
Without loss of generality we can assume that the hypercircle
tangent to the equator of a Bricard-Stachel octahedron $\mathcal{O}$
of type 3, passes through the (unique) point $\infty$ at infinity of
$\mathbb{H}^2$ which does not lie on the Euclidean line $z=0$, and
through the point $O=(0,0)$ at infinity of $\mathbb{H}^2$. Every
such hypercircle is specified by the equation $z=\rho\tang\alpha$
for some $\alpha\in(0,\frac{\pi}{2})\cup(\frac{\pi}{2},\pi)$. By the
symmetry of $\mathbb{H}^2$ with respect to the straight line
$\rho=0$, it is sufficient to consider the family of hypercircles
$K=\{z=\rho\tang\alpha|\alpha\in(0,\frac{\pi}{2})\}$. Let
$K_{\alpha}\in K$. We will construct all possible quadrangles
tangent to $K_{\alpha}$ such that none of their vertices belongs to
$K_{\alpha}$, i.e., all quadrangles that can serve as equators of
$\mathcal{O}$. Then we will verify the validity of the
theorem~\ref{sum_of_equator_edges_theorem} for such quadrangles.

Let us study the quadrangles based on the straight lines
${\Lambda}_{l}^{A_1}$, ${\Lambda}_{r}^{A_1}$, ${\Lambda}_{l}^{A_2}$,
 ${\Lambda}_{r}^{A_2}$, tangent to $K_{\alpha}$, which are
realised in $\mathbb{H}^2$ as the Euclidean demi-circles with the
centers $O^{A_1}_{l}=({\rho}_{A_1,l},0)$,
$O^{A_1}_{r}=({\rho}_{A_1,r},0)$, $O^{A_2}_{l}=({\rho}_{A_2,l},0)$,
$O^{A_2}_{r}=({\rho}_{A_2,r},0)$. Also, let ${\Lambda}_{l}^{A_1}$
and ${\Lambda}_{r}^{A_1}$ intersect at a point $A_1$,
${\Lambda}_{l}^{A_2}$ and ${\Lambda}_{r}^{A_2}$ intersect at a point
$A_2$. Assume that $A_1$ and $A_2$ are two opposite vertices of
$\mathcal{O}$, and that the inequalities
$0<{\rho}_{A_1,l}<{\rho}_{A_1,r}$, $0<{\rho}_{A_2,l}<{\rho}_{A_2,r}$
hold true.

\begin{remark}\label{r2} Let $T=({\rho}_{T}, z_{T})$ be a point in
$\mathbb{H}^2$, which serves as the intersection of straight lines
$\Lambda_{l}^{T}$ and $\Lambda_{r}^{T}$ tangent to a hypercircle
$K_{\alpha}$, and let $\Lambda_{l}^{T}$ and $\Lambda_{r}^{T}$ are
realised in $\mathbb{H}^2$ as the Euclidean demi-circles with the
centers  $O^{T}_{l}=({\rho}_{T,l},0)$, $O^{T}_{r}=({\rho}_{T,r},0)$
(${\rho}_{T,l}<{\rho}_{T,r}$). Then, by Remark~\ref{r1}, the angles
$\varphi_{T}^{l}\stackrel{\mathrm{def}}{=}\angle TO^{T}_{l}\rho$ and
$\phi_{T}^{r}\stackrel{\mathrm{def}}{=}\angle TO^{T}_{r}\rho$
determine uniquely the positions of $T$ on ${\Lambda}_{l}^{T}$ and
${\Lambda}_{r}^{T}$ correspondingly. Moreover,
\begin{equation}\label{cos_phi_T_l_and_r_hypercircle}
\cos\varphi_{T}^{l}=\frac{{\rho}_{T,r}}{{\rho}_{T,l}}\frac{{\cos}^{2}\alpha}{2\sin\alpha}
-\frac{1}{2\sin\alpha}-\frac{\sin\alpha}{2}\quad\mbox{and}\quad
\cos\varphi_{T}^{r}=\frac{{\rho}_{T,l}}{{\rho}_{T,r}}\frac{{\cos}^{2}\alpha}{2\sin\alpha}
-\frac{1}{2\sin\alpha}-\frac{\sin\alpha}{2}.
\end{equation}
\end{remark}

\textbf{Proof.} ${\Lambda}_{l}^{T}$ and ${\Lambda}_{r}^{T}$ are
tangent to $K_{\alpha}$. Hence, the radii $R_l$ and $R_r$ of the
demi-circles realizing ${\Lambda}_{l}^{T}$ and ${\Lambda}_{r}^{T}$
in $\mathbb{H}^2$ are determined by the formulas
\begin{equation}\label{R_l_and_R_r_hypercircle}
R_l={\rho}_{T,l}\sin\alpha\quad\mbox{and}\quad
R_r={\rho}_{T,r}\sin\alpha.
\end{equation}

Let $T_{\infty}$ be a point with coordinates $({\rho}_{T},0)$.
Applying the Euclidean Pythagorean theorem to $\triangle
TT_{\infty}O^{T}_{r}$ and simplifying the obtained expression, we
get:
\begin{equation}\label{Pythagoras_for_left triangle_3}
{\rho}_{T}^2+z_{T}^2=2{\rho}_{T}{\rho}_{T,l}-{\rho}_{T,l}^2{\cos}^{2}\alpha.
\end{equation}
Similarly, from $\triangle TT_{\infty}O^{T}_{l}$ we get that
\begin{equation}\label{Pythagoras_for_right triangle}
{\rho}_{T}^2+z_{T}^2=2{\rho}_{T}{\rho}_{T,r}-{\rho}_{T,r}^2{\cos}^{2}\alpha.
\end{equation}
Subtracting~(\ref{Pythagoras_for_left triangle_3})
from~(\ref{Pythagoras_for_right triangle}), we easily deduce:
\begin{equation}\label{rho_T_remarque_2_2}
{\rho}_{T}=\frac{{\rho}_{T,r}+{\rho}_{T,l}}{2}{\cos}^{2}\alpha.
\end{equation}

From the definitions of the cosines of $\varphi_{T}^{l}$ and
$\varphi_{T}^{r}$
($\cos\varphi_{T}^{l}=({\rho}_{T}-{\rho}_{T,l})/R_l$ and
$\cos\varphi_{T}^{r}=({\rho}_{T}-{\rho}_{T,r})/R_r$), taking into
account~(\ref{R_l_and_R_r_hypercircle})
and~(\ref{rho_T_remarque_2_2}), we
obtain~(\ref{cos_phi_T_l_and_r_hypercircle}). $\square$

By Remark~\ref{r2}, the angles
$\varphi_{A_1}^{{\Lambda}_{l}^{A_1}}\stackrel{\mathrm{def}}{=}\angle
A_1 O^{A_1}_{l}\rho$ and
$\varphi_{A_1}^{{\Lambda}_{r}^{A_1}}\stackrel{\mathrm{def}}{=}\angle
A_1 O^{A_1}_{r}\rho$ determine uniquely the positions of $A_1$ on
${\Lambda}_{l}^{A_1}$ and ${\Lambda}_{r}^{A_1}$ correspondingly.
Moreover,
\begin{equation}\label{cos_phi_A_1_l_and_r_hypercircle}
\cos\varphi_{A_1}^{{\Lambda}_{l}^{A_1}}=\frac{{\rho}_{A_1,r}}{{\rho}_{A_1,l}}\frac{{\cos}^{2}\alpha}{2\sin\alpha}
-\frac{1}{2\sin\alpha}-\frac{\sin\alpha}{2}\quad\mbox{and}\quad
\cos\varphi_{A_1}^{{\Lambda}_{r}^{A_1}}=\frac{{\rho}_{A_1,l}}{{\rho}_{A_1,r}}\frac{{\cos}^{2}\alpha}{2\sin\alpha}
-\frac{1}{2\sin\alpha}-\frac{\sin\alpha}{2}.
\end{equation}
Similarly, the angles
$\varphi_{A_2}^{{\Lambda}_{l}^{A_2}}\stackrel{\mathrm{def}}{=}\angle
A_2 O^{A_2}_{l}\rho$ and
$\varphi_{A_2}^{{\Lambda}_{r}^{A_2}}\stackrel{\mathrm{def}}{=}\angle
A_2 O^{A_2}_{r}\rho$ serve as the coordinates of $A_2$ on
${\Lambda}_{l}^{A_2}$ and ${\Lambda}_{r}^{A_2}$ correspondingly.
Moreover,
\begin{equation}\label{cos_phi_A_2_l_and_r_hypercircle}
\cos\varphi_{A_2}^{{\Lambda}_{l}^{A_2}}=\frac{{\rho}_{A_2,r}}{{\rho}_{A_2,l}}\frac{{\cos}^{2}\alpha}{2\sin\alpha}
-\frac{1}{2\sin\alpha}-\frac{\sin\alpha}{2}\quad\mbox{and}\quad
\cos\varphi_{A_2}^{{\Lambda}_{r}^{A_2}}=\frac{{\rho}_{A_2,l}}{{\rho}_{A_2,r}}\frac{{\cos}^{2}\alpha}{2\sin\alpha}
-\frac{1}{2\sin\alpha}-\frac{\sin\alpha}{2}.
\end{equation}

Suppose that the straight lines ${\Lambda}_{l}^{A_1}$ and
${\Lambda}_{l}^{A_2}$ intersect at a point $B_1$. Then the angles
$\varphi_{B_1}^{{\Lambda}_{l}^{A_1}}\stackrel{\mathrm{def}}{=}\angle
B_1 O^{A_1}_{l}\rho$ and
$\varphi_{B_1}^{{\Lambda}_{l}^{A_2}}\stackrel{\mathrm{def}}{=}\angle
B_1 O^{A_2}_{l} \rho$ serve as the coordinates of $B_1$ on
${\Lambda}_{l}^{A_1}$ and ${\Lambda}_{l}^{A_2}$ correspondingly.
Moreover,
\begin{equation}\label{cos_phi_B_1_l_and_r_hypercircle}
\cos\varphi_{B_1}^{{\Lambda}_{l}^{A_1}}=\frac{{\rho}_{A_2,l}}{{\rho}_{A_1,l}}\frac{{\cos}^{2}\alpha}{2\sin\alpha}
-\frac{1}{2\sin\alpha}-\frac{\sin\alpha}{2}\quad\mbox{and}\quad
\cos\varphi_{B_1}^{{\Lambda}_{l}^{A_2}}=\frac{{\rho}_{A_1,l}}{{\rho}_{A_2,l}}\frac{{\cos}^{2}\alpha}{2\sin\alpha}
-\frac{1}{2\sin\alpha}-\frac{\sin\alpha}{2}.
\end{equation}

Suppose also that ${\Lambda}_{r}^{A_1}$ and ${\Lambda}_{r}^{A_2}$
intersect at a point $B_2$. Then the angles
$\varphi_{B_2}^{{\Lambda}_{r}^{A_1}}\stackrel{\mathrm{def}}{=}\angle
B_2 O^{A_1}_{r}\rho$ and
$\varphi_{B_2}^{{\Lambda}_{r}^{A_2}}\stackrel{\mathrm{def}}{=}\angle
B_2 O^{A_2}_{r} \rho$ serve as the coordinates of $B_2$ on
${\Lambda}_{r}^{A_1}$ and ${\Lambda}_{r}^{A_2}$ correspondingly.
Moreover,
\begin{equation}\label{cos_phi_B_2_l_and_r_hypercircle}
\cos\varphi_{B_2}^{{\Lambda}_{r}^{A_1}}=\frac{{\rho}_{A_2,r}}{{\rho}_{A_1,r}}\frac{{\cos}^{2}\alpha}{2\sin\alpha}
-\frac{1}{2\sin\alpha}-\frac{\sin\alpha}{2}\quad\mbox{and}\quad
\cos\varphi_{B_2}^{{\Lambda}_{r}^{A_2}}=\frac{{\rho}_{A_1,r}}{{\rho}_{A_2,r}}\frac{{\cos}^{2}\alpha}{2\sin\alpha}
-\frac{1}{2\sin\alpha}-\frac{\sin\alpha}{2}.
\end{equation}

Suppose that ${\Lambda}_{r}^{A_1}$ and ${\Lambda}_{l}^{A_2}$
intersect at a point $C_1$. Then the angles
$\varphi_{C_1}^{{\Lambda}_{r}^{A_1}}\stackrel{\mathrm{def}}{=}\angle
C_1 O^{A_1}_{r}\rho$ and
$\varphi_{C_1}^{{\Lambda}_{l}^{A_2}}\stackrel{\mathrm{def}}{=}\angle
C_1 O^{A_2}_{l} \rho$ serve as the coordinates of $C_1$ on
${\Lambda}_{r}^{A_1}$ and ${\Lambda}_{l}^{A_2}$ correspondingly.
Moreover,
\begin{equation}\label{cos_phi_C_1_l_and_r_hypercircle}
\cos\varphi_{C_1}^{{\Lambda}_{l}^{A_2}}=\frac{{\rho}_{A_1,r}}{{\rho}_{A_2,l}}\frac{{\cos}^{2}\alpha}{2\sin\alpha}
-\frac{1}{2\sin\alpha}-\frac{\sin\alpha}{2}\quad\mbox{and}\quad
\cos\varphi_{C_1}^{{\Lambda}_{r}^{A_1}}=\frac{{\rho}_{A_2,l}}{{\rho}_{A_1,r}}\frac{{\cos}^{2}\alpha}{2\sin\alpha}
-\frac{1}{2\sin\alpha}-\frac{\sin\alpha}{2}.
\end{equation}

Suppose also that ${\Lambda}_{l}^{A_1}$ and ${\Lambda}_{r}^{A_2}$
intersect at a point $C_2$. Then the angles
$\varphi_{C_2}^{{\Lambda}_{l}^{A_1}}\stackrel{\mathrm{def}}{=}\angle
C_2 O^{A_1}_{l}\rho$ and
$\varphi_{C_2}^{{\Lambda}_{r}^{A_2}}\stackrel{\mathrm{def}}{=}\angle
C_2 O^{A_2}_{r} \rho$ serve as the coordinates of $C_2$ on
${\Lambda}_{l}^{A_1}$ and ${\Lambda}_{r}^{A_2}$ correspondingly.
Moreover,
\begin{equation}\label{cos_phi_C_2_l_and_r_hypercircle}
\cos\varphi_{C_2}^{{\Lambda}_{l}^{A_1}}=\frac{{\rho}_{A_2,r}}{{\rho}_{A_1,l}}\frac{{\cos}^{2}\alpha}{2\sin\alpha}
-\frac{1}{2\sin\alpha}-\frac{\sin\alpha}{2}\quad\mbox{and}\quad
\cos\varphi_{C_2}^{{\Lambda}_{r}^{A_2}}=\frac{{\rho}_{A_1,l}}{{\rho}_{A_2,r}}\frac{{\cos}^{2}\alpha}{2\sin\alpha}
-\frac{1}{2\sin\alpha}-\frac{\sin\alpha}{2}.
\end{equation}

As in the case of the quadrangles tangent to a horocycle
in~$\mathbb{H}^2$, the lengths of the edges of $A_1 B_1 A_2 B_2$ are
expressed
in~(\ref{length_of_A_1_B_1_horocycle})---(\ref{length_of_B_2_A_2_horocycle}),
and the lengths of the edges of $A_1 C_1 A_2 C_2$ are calculated
in~(\ref{length_of_C_1_A_1_horocycle})---(\ref{length_of_A_2_C_2_horocycle}).
Taking into
account~(\ref{cos_phi_A_1_l_and_r_hypercircle})---(\ref{cos_phi_C_2_l_and_r_hypercircle}),
it is easy to state the validity
of~(\ref{equator_A_1_B_1_A_2_B_2_horocycle})
and~(\ref{equator_A_1_C_1_A_2_C_2_horocycle}).

According to~(\ref{equator_A_1_B_1_A_2_B_2_horocycle})
and~(\ref{equator_A_1_C_1_A_2_C_2_horocycle}), the
theorem~\ref{sum_of_equator_edges_theorem} is valid for any equator
of a Bricard-Stachel octahedron of type 3 tangent to a hypercircle
in at least one of its flat positions.

\bigskip

The case when three vertices of an equator of a flexible octahedron
in its flat position lie on a straight line, is similar. The case
when all four vertices of an equator lie on a straight line, is
trivial.

\bigskip

The author is grateful to Victor Aleksandrov for his help at all
phases of preparation of this paper.

Dmitriy Slutskiy

Sobolev Institute of Mathematics of the SB RAS,

4 Acad. Koptyug avenue, 630090 Novosibirsk, Russia

and

Novosibirsk State University,

2 Pirogova Street, 630090, Novosibirsk, Russia

slutski@ngs.ru

\end{document}